\documentclass[11pt,a4paper,dvipsnames]{article}

\usepackage[top=25mm,bottom=25mm,left=25mm,right=25mm]{geometry}

\usepackage[T1]{fontenc}
\usepackage[utf8]{inputenc}
\usepackage[english]{babel}
\usepackage{lmodern}
\usepackage{microtype}

\usepackage{amsmath,amssymb,amsfonts,mathtools,amsthm}
\usepackage{bm}
\usepackage{arydshln}

\usepackage[inline,shortlabels]{enumitem}

\usepackage[dvipsnames]{xcolor}

\definecolor{linkblue}{RGB}{0,0,220}
\definecolor{tocblue}{RGB}{0,0,220}
\definecolor{referenceblue}{RGB}{0,0,220}

\colorlet{MidnightBlue}{linkblue}

\usepackage{graphicx}
\usepackage{float}
\usepackage{tikz}
\usepackage{pgfplots}
\pgfplotsset{compat=1.18}

\usepackage{booktabs,tabularx,array}

\usepackage[bf,sf]{titlesec}
\usepackage[labelfont={bf,sf}]{caption}
\usepackage{subcaption}
\usepackage{csquotes}

\usepackage[
backend=biber,
style=ieee,
maxnames=20,
giveninits=true,
isbn=false,
doi=true,
url=true
]{biblatex}

\usepackage{hyperref}

\hypersetup{
	colorlinks=true,
	citecolor=referenceblue,
	linkcolor=tocblue,
	urlcolor=linkblue,
	hypertexnames=true,
	plainpages=false
}

\usepackage[capitalize,nameinlink]{cleveref}

\numberwithin{equation}{section}

\makeatletter

\renewcommand{\tagform@}[1]{%
	\maketag@@@{%
		{\color{tocblue}%
			(\ignorespaces#1\unskip\@@italiccorr)}%
	}%
}

\renewcommand{\maketitle}{%
	\begin{center}
		\vspace*{1.2cm}
		{\sffamily\bfseries
			\fontsize{16.5}{20}\selectfont
			\@title\par}
		\vspace{10pt}
		{\large\@author\par}
		\vspace{19pt}
		{\normalsize\@date\par}
	\end{center}
	\vspace{8pt}
}

\makeatother

\providecommand{\keywords}[1]{%
	\par\noindent
	\textbf{\textsf{Keywords.}} #1\par
}

\theoremstyle{plain}

\newtheorem{theorem}{Theorem}
\newtheorem{lemma}{Lemma}

\theoremstyle{definition}

\newtheorem{assumption}{Assumption}
\newtheorem{remark}{Remark}
\newtheorem{problem}{Problem}
\newtheorem{example}{Example}

\newcommand{\qedblack}{\hfill$\blacksquare$}

\newcommand{\overbar}[1]{%
	\mkern 1.5mu
	\overline{\mkern-1.5mu#1\mkern-1.5mu}%
	\mkern 1.5mu
}

\renewcommand{\arraystretch}{0.8}
\title{Safe and Scalable Minimax Adaptive Control of Positive Systems}
\author{Fethi Bencherki$^{\star}$\qquad David Ohlin$^{\star}$\qquad Tomas Meijer$^{\star}$}
\date{$^{\star}$Department of Automatic Control, Lund University, Lund, Sweden\\[2pt]
{\color{MidnightBlue}\ttfamily\small\{fethi.bencherki,david.ohlin,tomas.meijer\}@control.lth.se}}

\begin{document}
\maketitle

\begin{abstract}
	This paper studies minimax adaptive control of positive linear systems
	with fixed unknown dynamics and nonnegative disturbances.
	We reformulate the problem by compressing the historical data into
	a vector and a positivity-induced consistency set containing the
	models compatible with the observations.
	Together with the physical state, these quantities form a sufficient
	information state and yield an exact reformulation as a standard minimax
	dynamic game.
	We allow the admissible input set to adapt to the current
	consistency set, ensuring positivity for every
	remaining feasible model. As inconsistent models are eliminated,
	the input set can expand, allowing less cautious control actions.
	For finite model sets, we derive an explicit solution to the
	Bellman inequality and construct a safe and stabilizing
	certainty-equivalence controller with robustness guarantees.
	The framework also accommodates structured classes of continuum
	model uncertainty whose consistency sets admit fixed,
	finite-dimensional representations, enabling scalable
	information-state updates.
\end{abstract}

\section{Introduction}
Modern controllers are increasingly deployed in systems
with uncertain, poorly modeled, or changing dynamics, making first-principles modeling difficult. System identification
and data-driven control use measured data to sidestep this difficulty, but
typically separate data collection from control design
\cite{Ljung1999,vanWaarde2023}. Adaptive control instead updates its
decisions as observations become available~\cite{AAstrom2013}.
This interaction creates the dual control problem: an input affects both
current performance and what can be learned about the system. Dynamic
programming provides a natural framework to address such problems, but its information state
typically contains the complete observation history or a complicated belief
state, making solutions difficult to obtain
\cite{Feldbaum1960,Wittenmark1995,Mesbah2018,meijer2026dual}.

Minimax adaptive control provides a deterministic alternative to stochastic
formulations of dual control. It extends the worst-case perspective of
$\mathcal H_\infty$ control \cite{zames1981feedback,zhou1996robust,bacsar2008h}
to settings in which the dynamics must be learned online, treating both
disturbances and unknown parameters as adversarial quantities
\cite{vinnicombe2004examples,megretski2004nonlinear}.
These nonlinear adaptive controllers can stabilize systems
with uncertainties beyond those that can be handled
by robust linear time-invariant (LTI) controllers
\cite{bencherki2023robust}.
An early formulation for linear systems appeared in
\cite{didinsky1994minimax}, followed by extensions to nonlinear systems
in \cite{pan1998adaptive}. A central challenge is to replace the growing
measurement history with a recursively updated, finite-dimensional
information state suitable for Bellman recursion. For a finite family
of linear models, \cite{rantzer2021minimax} showed that a quadratic
data covariance matrix provides such a representation: augmenting the physical
state with this matrix transforms the adaptive control problem into
a two-player zero-sum dynamic game that can be analyzed using
minimax dynamic programming.

Positive systems provide a particularly compelling setting for dual control
because the nonnegativity constraints limit how the system can be probed, while
the resulting order structure offers opportunities for tractable analysis.
They describe flow, transportation, biological, and compartmental processes
in which the states, and often the inputs, represent nonnegative quantities.
An input may add material, energy, or population whose effect cannot be
immediately reversed, so learning must be balanced carefully against current
performance and safety. At the same time, positivity can simplify worst-case
comparisons and enable scalable control design~\cite{Rantzer2021b}. These
aspects have motivated recent work on model-based optimal and minimax control
~\cite{Gurpegui2023}, as well as data-driven and
adaptive methods~\cite{miller2023data,Iwata2025,bencherki2025adaptive,bencherki2024data}.

The formulations in~\cite{Bencherki2026,Bencherki2026b} consider
positivity-preserving disturbances. Under this disturbance model, every
candidate system can explain a nonnegative next-state observation by selecting
an admissible disturbance. Consequently, candidate models cannot be eliminated
by testing their consistency with the observed data. Instead, their accumulated
disturbance costs are tracked using model-specific history variables. This construction maintains one history variable per candidate model, so the dimension of the resulting representation depends on the size of the model family. In the linear quadratic setting, previous work on minimax adaptive
dual control \cite{rantzer2021minimax} uses a covariance-based
information state to compress the required memory to a fixed
finite dimension. In the present setting, we prove that the
historical information required for control
can be represented by a single historical data vector and a
positivity-induced consistency set.
Related work in the same positive system setting has also shown that,
when the system model is known, a vector whose dimension matches the
system's state dimension can suffice for state estimation and
control~\cite{ohlin26thesis}.

\textbf{Contributions.}
The present paper considers nonnegative disturbances, which make observed
transitions informative: a model is eliminated whenever it cannot explain an
observation using a nonnegative disturbance. We use this observation to
replace model-specific histories with a consistency set and an accumulated
historical data vector. The main contributions are as follows:

 \begin{itemize}
\item[\textit{(C1)}] We formulate a minimax adaptive control problem for
positive systems in which the adversary selects fixed unknown dynamics and
nonnegative disturbances, while the controller selects the inputs; see
Section~\ref{sec:problem-formulation}.

\item[\textit{(C2)}] We derive an exact reformulation as a standard
two-player minimax dynamic game and show that the physical state,
consistency set, and historical-data vector form a sufficient information
state for dynamic programming; see Sections~\ref{sec:reformulation}
and~\ref{sec:minimax-dp}.

\item[\textit{(C3)}] For finite model sets with model-dependent input bounds,
we construct a safe CE controller and an explicit solution to the Bellman
inequality. The consistency set yields cautious control actions initially,
ensuring positivity for every model remaining in the uncertainty set, while
permitting less cautious control as uncertainty shrinks; see
Section~\ref{sec:finite-model-certificate}.

\item[\textit{(C4)}]
We establish a scalable information-state representation
for structured continuum uncertainty classes whose
consistency sets admit fixed, finite-dimensional recursive
representations. Together with a common historical-data
vector, this yields an information state whose dimension
is independent of the number of candidate model pairs,
enabling efficient online updates even when the
uncertainty set is uncountable; see
Section~\ref{sec:history-compression}.

\item[\textit{(C5)}] We illustrate and validate the theoretical results
through numerical experiments; see
Section~\ref{sec:numerical-experiments}.
\end{itemize}
\textbf{Organization.}
Section~\ref{sec:problem-formulation} presents the problem formulation.
Section~\ref{sec:reformulation} derives a data-history reformulation that
converts the original problem into a standard minimax dynamic game.
Section~\ref{sec:minimax-dp} analyzes this game using minimax dynamic
programming. Section~\ref{sec:finite-model-certificate} presents explicit
solutions to the Bellman inequality for a finite set of positive LTI plants.
Section~\ref{sec:history-compression} identifies continuum model classes for
which the hyperstate admits a fixed-dimensional recursive representation.
Section~\ref{sec:numerical-experiments} illustrates the theoretical results
through numerical examples. Finally, Section~\ref{sec:conclusion} concludes
the paper and discusses directions for future research. Proofs of the main
results are deferred to the appendix.

\noindent\textbf{Notation.}
 Let \(\mathbb R_+=[0,\infty)\) and \(\mathbb N=\{0,1,2,\ldots\}\). Vector and matrix inequalities and absolute values are interpreted elementwise. For a matrix \(M\), \([M]_{q\ell}\) denotes its \((q,\ell)\)th entry, while \([a]_q\) denotes the \(q\)th entry of a vector \(a\). The maximum and minimum of two vectors are taken elementwise; that is,
$
[\max\{a,b\}]_q=\max\{[a]_q,[b]_q\},
$
and similarly for \(\min\{a,b\}\). For a vector \(a\), \(\operatorname{diag}(a)\) denotes the diagonal matrix with diagonal entries given by \(a\), and \(\operatorname{sgn}(a)\) denotes its elementwise sign, with \(\operatorname{sgn}(0)=0\). The symbols \(I\), \(\mathbf 0\), and \(\mathbf 1\) denote the identity matrix, zero vector or matrix, and all-ones vector or matrix of appropriate dimensions, respectively. The transpose of a matrix \(M\) is denoted by \(M^\top\), and \(\rho(M)\) denotes its spectral radius. For a finite set \(\mathcal C\), \(|\mathcal C|\) denotes its cardinality.
\section{Problem formulation}
\label{sec:problem-formulation}

\subsection{Positive system class}

We consider the uncertain discrete-time positive system
\begin{equation}
    x_{t+1}=Ax_t+Bu_t+w_t,
    \qquad x_0\in\mathbb R_+^n,
    \label{eq:system}
\end{equation}
where the unknown pair $(A_\star,B_\star)$ belongs to a known compact set
$\mathcal M\subset\mathbb R_+^{n\times n}\times\mathbb R^{n\times m}$ as stated in the next assumption.
\begin{assumption}[Model uncertainty]
\label{ass:model-uncertainty}
The model set
$\mathcal M\subset\mathbb R_+^{n\times n}\times\mathbb R^{n\times m}$
is known and compact, and the true system matrices satisfy
$(A_\star,B_\star)\in\mathcal M$.
\end{assumption}
Unlike the positivity-preserving disturbance class
in~\cite{Bencherki2026}, disturbances here must be nonnegative:
\begin{equation}
    w_t\geq0,
    \qquad t\in\mathbb N.
    \label{eq:nonnegative-disturbance}
\end{equation}
Given the observed state--input history
$x_0,\ldots,x_t$ and $u_0,\ldots,u_{t-1}$, define
\begin{align}\label{const-set}
\mathcal C_t \coloneqq
\bigl\{(A,B)\in\mathcal M:\;
x_{\tau+1}\geq Ax_\tau+Bu_\tau,
0\leq\tau<t\bigr\},
\qquad \mathcal C_0=\mathcal M.
\end{align}
Thus, $
\mathcal M\supseteq\mathcal C_t
$ is the set of fixed models compatible with all
transitions observed before time $t$. For every nonempty consistency set,
we define the set of admissible inputs $\mathcal U(x,\mathcal C)$ by
\begin{equation}
\mathcal U(x,\mathcal C)\subseteq
\bigcap_{(A,B)\in\mathcal C}
\{u\in\mathbb R^m:Ax+Bu\geq0\}.
\label{eq:input-set}
\end{equation}
It follows from \eqref{eq:input-set}--\eqref{eq:nonnegative-disturbance} that
$x_{t+1}\geq0$ whenever $x_t\geq0$.
Thus, every admissible input preserves positivity for all data-consistent
models, including the true pair since
$(A_\star,B_\star)\in\mathcal C_t$ for all $t \in \mathbb N$. For example, let
$\mathcal M=\{(A_1,B_1),(A_2,B_2)\}$. While both models remain consistent,
the input must preserve positivity for both; once model $i$ is identified,
only its positivity constraint remains. The admissible input set may
therefore expand as the consistency set shrinks, permitting less cautious
control actions. The model-based counterpart of the problem considered here has been addressed in~\cite{Gurpegui2023}.
\subsection{A minimax dual control problem}
The controller observes the state history and its past inputs, but neither
the true dynamics nor the disturbances, and applies a deterministic causal
policy
\begin{equation}
u_t=\mu_t(x_0,\ldots,x_t;u_0,\ldots,u_{t-1}) \in\mathcal U(x_t,\mathcal C_t).
\label{eq:causal-policy}
\end{equation}
For a given $s\in\mathbb R_+^n$, $r\in\mathbb R^m$, and
$\gamma\in\mathbb R_+^n$, in this paper, we consider the following minimax dynamic game problem.
\begin{problem}
   Find a causal policy $\mu$ solving
\begin{equation}
\begin{aligned}
\!\! \!\! J_\star(x_0)
\coloneqq{}&
\inf_{\mu}
\sup_{\substack{(A,B)\in\mathcal M, w\geq0\\T\in\mathbb N}}
\sum_{t=0}^{T}
\left(
s^\top x_t+r^\top u_t-\gamma^\top w_t
\right) \\
\mathrm{s.t.}\quad
&x_{t+1}=Ax_t+Bu_t+w_t,\qquad x_0\in\mathbb R_+^n,\\
&u_t=\mu_t(x_0,\ldots,x_t;u_0,\ldots,u_{t-1})
\in\mathcal U(x_t,\mathcal C_t).
\end{aligned}
\label{eq:problem-1}
\end{equation}
\label{pb1}
\end{problem}
We impose the following assumption on the stage cost.
\begin{assumption}[Positivity of the stage cost]
\label{ass:positive-stage-cost}
Let $s,\gamma\in\mathbb R_+^n$ and $r\in\mathbb R^m$. For every
$x\in\mathbb R_+^n$ and
$u\in\mathcal U(x,\mathcal C)$, assume
$
\max_{w\geq0}\{s^\top x+r^\top u-\gamma^\top w\}\geq0.
\label{eq:positive-stage-cost}
$
\end{assumption}
\begin{remark}
Since $\gamma,w\in\mathbb R_+^n$, the term $-\gamma^\top w$ is maximized at
$w=0$. Consequently, Assumption~\ref{ass:positive-stage-cost} is equivalent to
$
s^\top x+r^\top u\geq0.
$
\end{remark}
\section{Reformulated Minimax Problem}
\label{sec:reformulation}
We reformulate Problem~\ref{pb1} as a standard minimax dynamic game by
treating the next state as the adversary's action and augmenting the state
with the consistency set and historical-data vector.
\subsection{Reparameterization through the next state}
Since $w_t\geq0$, the disturbance can steer the next state to any
$v_t\geq Ax_t+Bu_t$. We may therefore replace $w_t$ with the adversarial
next-state variable
\[
v_t\coloneqq x_{t+1},
\qquad
v_t\geq Ax_t+Bu_t.
\]
Applying this change of variables to Problem~\ref{pb1} gives
\begin{equation}
\begin{aligned}
J_\star(x_0)
\!\coloneqq{}&
\!\!\inf_{\mu}\!\!\!
\sup_{\substack{
(A,B)\in\mathcal M\\
v\geq Ax+Bu\\
T\in\mathbb N}}
\!\sum_{t=0}^{T}
s^\top x_t+r^\top u_t
\!-\gamma^\top(v_t-Ax_t-Bu_t)
\\
\mathrm{s.t.}\quad
&x_{t+1}=v_t,
\qquad x_0\in\mathbb R_+^n,
\\
&u_t=\!\mu_t(x_0,\ldots,x_t;\!u_0,\ldots,u_{t-1})
\!\in\!\mathcal U(x_t,\mathcal C_t).
\end{aligned}
\label{eq:problem-after-v}
\end{equation}
The state equation in~\eqref{eq:problem-after-v} no longer contains $(A,B)$
explicitly, but the set of feasible next-state trajectories remains
model-dependent. Replacing the constraint
$v_t\geq Ax_t+Bu_t$ with $v_t\geq0$ would enlarge the adversary's feasible
set. Consequently, unlike in~\cite{Bencherki2026,Bencherki2026b}, the maximizations over the
trajectory and the model cannot be separated. Instead, they can be expressed
in terms of the models consistent with each trajectory, as shown next.
\subsection{Positivity-induced consistency set and history vector}
From~\eqref{const-set}, the consistency set evolves recursively according to
\begin{equation}
\mathcal C_{t+1}
\coloneqq
\left\{
(A,B)\in\mathcal C_t:
v_t\geq Ax_t+Bu_t
\right\}\subseteq \mathcal C_t,\quad \mathcal C_0\coloneqq\mathcal M
\label{eq:C-update}
\end{equation}
Equation~\eqref{eq:C-update} implies
\begin{equation}
(A,B)\in\mathcal C_{t+1}
\iff (A,B)\in\mathcal M,\
v_\tau\geq Ax_\tau+Bu_\tau,\ 0\leq\tau\leq t.
\label{eq:C-meaning}
\end{equation}
Thus, every element of $\mathcal C_{t+1}$ is one fixed model that is
consistent with the entire observed trajectory up to time $t$.
Next, define the historical-data vector $Z_t \in \mathbb R^{2n+m}$ and its recursion by
\begin{equation}
Z_t\coloneqq
\sum_{\tau=0}^{t-1}
\begin{bmatrix}
v_\tau\\
x_\tau\\
u_\tau
\end{bmatrix}
\quad\Longrightarrow\quad
Z_{t+1}=Z_t+
\begin{bmatrix}
v_t\\
x_t\\
u_t
\end{bmatrix},
\qquad
Z_0=0.
\label{eq:Z-definition}
\end{equation}
For every $(A,B)\in\mathcal C_{T+1}$, the disturbance sequence associated with
that model satisfies
$w_t=v_t-Ax_t-Bu_t\geq0$ for $0\leq t\leq T$. Consequently,
\begin{equation*}
\sum_{t=0}^{T}\gamma^\top w_t
=\sum_{t=0}^{T}\gamma^\top(v_t-Ax_t-Bu_t)
=\gamma^\top[\,I\ {-}A\ {-}B\,]Z_{T+1}.
\end{equation*}
Thus, the accumulated disturbance penalty associated with every
consistent model can be recovered from the same vector $Z_{T+1}$.
For a fixed policy and horizon, the model and trajectory maximization in
\eqref{eq:problem-after-v} is
\begin{equation}
\begin{aligned}
\sup_{\substack{
(A,B)\in\mathcal M\\
v\geq Ax+Bu
}}
\sum_{t=0}^{T}
\left(s^\top x_t+r^\top u_t\right)
-
\gamma^\top
\begin{bmatrix}
I & -A & -B
\end{bmatrix}
Z_{T+1}.
\end{aligned}
\label{eq:before-regrouping}
\end{equation}
Taking \(t=T\) in \eqref{eq:C-meaning} gives
\begin{equation}
\mathcal C_{T+1}\!=\bigl\{(A,B)\in\mathcal M:
v_t\geq Ax_t+Bu_t\ \forall\ 0\leq t\leq T\bigr\}.
\label{eq:C-final-meaning}
\end{equation}
The following result shows that the joint model--trajectory supremum in
\eqref{eq:before-regrouping} can be expressed equivalently using the final
consistency set, without changing the adversary's feasible action set.
\begin{lemma}[Equivalence of feasible descriptions]\label{lemma1}
Fix a policy and a horizon \(T\). A model--trajectory pair \(\bigl((A,B),(v_0,\ldots,v_T)\bigr)\) is feasible in \eqref{eq:before-regrouping}, meaning that \((A,B)\in\mathcal M\) and \(v_\tau\geq Ax_\tau+Bu_\tau\) for every \(\tau=0,\ldots,T\), if and only if \(\mathcal C_{T+1}\neq\varnothing\) and \((A,B)\in\mathcal C_{T+1}\). Thus, the original and regrouped descriptions contain exactly the same feasible model--trajectory pairs.
\end{lemma}

{\em Proof.}
Suppose first that the model--trajectory pair is feasible in \eqref{eq:before-regrouping}. Then \((A,B)\in\mathcal M\) and \(v_t\geq Ax_t+Bu_t\) for \(0\leq t\leq T\). By \eqref{eq:C-final-meaning}, it follows that \((A,B)\in\mathcal C_{T+1}\), and hence \(\mathcal C_{T+1}\neq\varnothing\). Conversely, suppose that \(\mathcal C_{T+1}\neq\varnothing\) and \((A,B)\in\mathcal C_{T+1}\). Since \(\mathcal C_{T+1}\subseteq\mathcal M\), we have \((A,B)\in\mathcal M\), while \eqref{eq:C-final-meaning} gives \(v_t\geq Ax_t+Bu_t\) for \(0\leq t\leq T\). Thus, the same model--trajectory pair satisfies all the constraints in \eqref{eq:before-regrouping} and is feasible.
\qedblack
\paragraph{An equivalent minimax dynamic game}

It follows from Lemma~\ref{lemma1} that the joint supremum in
\eqref{eq:before-regrouping} can be indexed by trajectories with a nonempty
final consistency set and, for each such trajectory, by the models belonging
to that consistency set. Furthermore, the consistency
 set dynamics and the
constraint
$
\mathcal C_{T+1}\neq\varnothing
$
incorporate the coupling between the fixed model and the trajectory. In particular, \(\mathcal C_{T+1}\) contains exactly the fixed models
that are consistent with the entire trajectory. Since the stage cost is independent of the model parameters
\((A,B)\in\mathcal C_{T+1}\), \eqref{eq:before-regrouping} is equivalent to
\begin{equation*}
\sup_{\substack{
v\\
\mathcal C_{T+1}\neq\varnothing
}} 
\biggl\{
\sum_{t=0}^{T}
\left(s^\top x_t+r^\top u_t\right)
- 
\min_{(A,B)\in\mathcal C_{T+1}} 
\gamma^\top
\begin{bmatrix}
I & -A & -B
\end{bmatrix}
Z_{T+1}\biggr\}.
\label{eq:regrouped-cost}
\end{equation*}
Thus, the fixed model enters the objective explicitly only through the
terminal penalty, while its compatibility with the trajectory is enforced
by the consistency set dynamics. Using \eqref{eq:C-update} together with the
hyperstate dynamics, Problem~\ref{pb1} can therefore be written as the
equivalent minimax dynamic game stated next.
\begin{problem}
Find a causal policy $\eta$ solving
\begin{equation}
	\begin{aligned}[c]
		J_\star(x_0)
		={}&
		\inf_{\eta}
		\sup_{\substack{
				v\\
				\mathcal C_{T+1}\neq\varnothing\\
				T\in\mathbb N}}
		\biggl\{
		\sum_{t=0}^{T}
		\left(s^\top x_t+r^\top u_t\right)
		-
		\min_{(A,B)\in\mathcal C_{T+1}}
		\gamma^\top[\,I\ {-}A\ {-}B\,]Z_{T+1}
		\biggr\}
		\\
		\mathrm{s.t.}\quad
		&x_{t+1}=v_t,
		\\
		&\mathcal C_{t+1}
		=
		\left\{
		(A,B)\in\mathcal C_t:
		v_t\geq Ax_t+Bu_t
		\right\},
		\quad
		\mathcal C_0=\mathcal M,
		\\
		&Z_{t+1}
		=
		Z_t+
		\begin{bmatrix}
			v_t^\top &
			x_t^\top &
			u_t^\top
		\end{bmatrix}^\top,
		\quad
		Z_0=0,
		\\
		&u_t
		=
		\eta(x_t,\mathcal C_t,Z_t),
		\quad
		u_t\in\mathcal U(x_t,\mathcal C_t).
	\end{aligned}
	\raisetag{4\baselineskip}
	\label{eq:problem-2}
\end{equation}
\label{pb2}
\end{problem}
The hyperstate of the reformulated game is
$(x_t,\mathcal C_t,Z_t)$. Given this hyperstate and the actions
$u_t$ and $v_t$, the recursions in \eqref{eq:problem-2} uniquely
determine its successor.
The update requires no knowledge of the true
pair. At each stage, the controller selects $u_t$, after which the
adversary selects $v_t$ subject to
$\mathcal C_{t+1}\neq\varnothing$. In general, however, the information state $\mathcal C_t$ is
infinite-dimensional and need not admit a finite-dimensional
parametrization. Section~\ref{sec:history-compression} identifies model
classes for which $\mathcal C_t$ admits a finite-dimensional recursive
representation.
Problem~\ref{pb2} is a deterministic
two-player minimax dynamic game. It will be addressed next via minimax dynamic programming.
\section{Minimax dynamic programming}
\label{sec:minimax-dp}

We now solve the reformulated problem~\ref{pb2} using minimax dynamic programming on
the hyperstate $(x,\mathcal C,Z)$.

\textbf{Dynamic programming operators.}
For a value function $V(x,\mathcal C,Z)$, define the operators
$\mathcal F$ and $\mathcal F_u$ by
\begin{equation}
	\mathcal F V(x,\mathcal C,Z)
	\coloneqq
	\min_{u\in\mathcal U(x,\mathcal C)}
	\underbrace{
		\max_{\substack{v,\ \mathcal C_+\neq\varnothing}}
		\left\{
		s^\top x+r^\top u+V(v,\mathcal C_+,Z_+)
		\right\}
	}_{\eqqcolon\mathcal F_uV(x,\mathcal C,Z)}.
	\label{eq:Bellman-operator}
\end{equation}
where $\mathcal C_+$ and $Z_+$ denote the updated consistency set and
history data vector, respectively, with the update
\begin{equation}
\left\{
\begin{aligned}
x_+&=v,\\
\mathcal C_+
&=\{(A,B)\in\mathcal C:v\geq Ax+Bu\},\\
Z_+&=Z+\begin{bmatrix}v^\top & x^\top &u^\top \end{bmatrix}^\top.
\end{aligned}
\right.
\label{eq:dp-hyperstate-update}
\end{equation}
The condition $\mathcal C_+\neq\varnothing$ means that the adversary may choose only
a next state that is compatible with at least one currently consistent
model.

\textbf{Initialization and value iteration.}
We initialize value iteration with the terminal cost function
\begin{equation}
V_0(x,\mathcal C,Z)
=
-\min_{(A,B)\in\mathcal C}
\gamma^\top[\,I\;-A\;-B\,]Z.
\label{eq:V0}
\end{equation}
The value iteration is then carried recursively by
\begin{align}
V_{k+1}(x,\mathcal C,Z)
=\mathcal F V_k(x,\mathcal C,Z)
=
\min_{u\in\mathcal U(x,\mathcal C)}
\max_{\substack{v, \mathcal C_+\neq\varnothing}}
\left\{
s^\top x+r^\top u+
V_k(v,\mathcal C_+,Z_+)
\right\}.
\label{eq:value-iteration}
\end{align}
This defines the sequence of value functions
$\{V_k\}_{k\in\mathbb N}$. The following result establishes the equivalence
between the original Problem~\ref{pb1} and its reformulation as a standard dynamic
game, Problem~\ref{pb2}.

\begin{theorem}
\label{thm:minimax-dp}
Let $s\in\mathbb R_+^n$, $r\in\mathbb R^m$, and
$\gamma\in\mathbb R_+^n$, and let $\mathcal M$ be a parameter
set such that Assumptions~\ref{ass:model-uncertainty}
and~\ref{ass:positive-stage-cost} hold. Define
$\mathcal F$ as in \eqref{eq:Bellman-operator} and the sequence
$\{V_k\}_{k\in\mathbb N}$ as in
\eqref{eq:V0}--\eqref{eq:value-iteration}. Then, the following statements hold:
\begin{enumerate}
\item[(i)]
Problems~\ref{pb1} and~\ref{pb2} have the same value. This value is finite if and only if
the sequence $\{V_k(x,\mathcal M,0)\}_{k=0}^{\infty}$ is bounded above
for every $x\in\mathbb R_+^n$. In that case, the limit
$V_\star\coloneqq\lim_{k\to\infty}V_k$ exists and satisfies
$J_\star(x_0)=V_\star(x_0,\mathcal M,0)$.
\item[(ii)]
Let $
\eta^\star(x,\mathcal C,Z)\in
\arg\min_{u\in\mathcal U(x,\mathcal C)}
\mathcal F_uV_\star(x,\mathcal C,Z).$
Then $\eta^\star$ is an optimal control law for Problem~\ref{pb2}, and the control
law
\begin{equation}
\mu_t^\star
(x_0,\ldots,x_t;u_0,\ldots,u_{t-1})
\coloneqq
\eta^\star(x_t,\mathcal C_t,Z_t),
\label{eq:mu-star}
\end{equation}
where
\begin{equation}
\begin{aligned}
\mathcal C_t
&=
\left\{
(A,B)\in\mathcal M:
x_{\tau+1}\geq Ax_\tau+Bu_\tau,\;
0\leq\tau<t
\right\},\\
Z_t
&=
\sum_{\tau=0}^{t-1}
\begin{bmatrix}
x_{\tau+1}^\top &
x_\tau ^\top &
u_\tau^\top
\end{bmatrix}^\top,
\end{aligned}
\label{eq:Ct-history}
\end{equation}
 is optimal for
Problem~\ref{pb1}.

\item[(iii)]
Suppose there exist an admissible feedback law
\(\overline\eta(x,\mathcal C,Z)\in\mathcal U(x,\mathcal C)\) and a function
\(\overline V(x,\mathcal C,Z)\) such that
\(\overline V\geq V_0\) and
\(\mathcal F_{\overline\eta(x,\mathcal{C},Z)}\overline V\leq\overline V\) for every
admissible \((x,\mathcal C,Z)\), where
\[
\mathcal F_{\overline\eta(x,\mathcal{C},Z)}V(x,\mathcal C,Z)
\coloneqq
\left.\mathcal F_uV(x,\mathcal C,Z)
\right|_{u=\overline\eta(x,\mathcal C,Z)}.
\]
Then, the policy \(\overline\mu\) defined by
\begin{equation}
u_t
\!=\!
\overline\mu_t(x_0,\ldots,x_t;u_0,\ldots,u_{t-1})
\coloneqq
\overline\eta(x_t,\mathcal C_t,Z_t),
\label{eq:policy-induced-by-hyperstate}
\end{equation}
where \(\mathcal C_t\) and \(Z_t\) are recursively constructed from~\eqref{eq:Ct-history}, satisfies
\(
J_{\overline\mu}(x_0)\leq
\overline V(x_0,\mathcal M,0).
\)
\end{enumerate}
\end{theorem}
{\em Proof.} See Appendix~\ref{proofappenA}. \qedblack

\begin{remark}
Theorem~\ref{thm:minimax-dp} has three main implications:
\begin{itemize}
    \item[\textit{(i)}] The hyperstate
$(x,\mathcal C,Z)$ contains all the information required for optimal
control. A policy obtained from the reformulated dynamic game induces a
causal policy for the original uncertain system because $\mathcal C_t$
and $Z_t$ can be updated recursively from measured data.
    \item[\textit{(ii)}] Any feasible solution of the associated Bellman equation/inequality for the reformulated problem~\ref{pb2} provides a
certified performance bound for the induced policy for the original problem~\ref{pb1}.
\item[\textit{(iii)}] 
In contrast to \cite{Bencherki2026,Bencherki2026b},
which consider positivity-preserving disturbances and employ
model-specific history variables, the present work restricts
the disturbance class to nonnegative disturbances.
Under this assumption, the information state can be
represented using a common historical-data vector and
a consistency set. As shown in
Section~\ref{sec:history-compression}, for certain continuum
uncertainty classes, the consistency set $\mathcal C_t$
admits a fixed, finite-dimensional recursive representation.
Consequently, the hyperstate remains finite-dimensional
even when $\mathcal M$ is uncountable.
\end{itemize}
\end{remark}

\begin{remark}[Online implementation]
The controller evaluates the feedback law at
$(x_t,\mathcal C_t,Z_t)$ and applies $u_t$. At time $t+1$, after observing
$x_{t+1}$, it updates
\begin{equation}
\begin{aligned}
\mathcal C_{t+1}
&=\{(A,B)\in\mathcal C_t:
x_{t+1}\geq Ax_t+Bu_t\},\\
Z_{t+1}
&=Z_t+
\begin{bmatrix}
x_{t+1}^\top& x_t^\top& u_t^\top
\end{bmatrix}^\top.
\end{aligned}
\label{eq:online-update}
\end{equation}
Both updates use only measured data and require no knowledge of the
true pair $(A_\star,B_\star)$, which remains in every consistency set since
$w_t=x_{t+1}-( A_\star x_t+B_\star u_t)\geq0$. The update simply discards
inconsistent candidate dynamics $(A,B)$.
\end{remark}
Having established the equivalence between Problems~\ref{pb1}
and~\ref{pb2}, we now construct a certainty-equivalence (CE)
controller for a finite set of candidate models and certify
its performance.
\section{A solution to the Bellman inequality for finite model sets}
\label{sec:finite-model-certificate}
Consider the case when the true system $(A_\star,B_\star)$ belongs to the finite model set
\begin{align*}
	\mathcal M=\{(A_1,B_1),\ldots,(A_M,B_M)\}.
\end{align*}
For each candidate model $(A_i,B_i)\in\mathcal M$, we make the following
assumption.
\begin{assumption}[Model-dependent input bounds]
\label{ass:model-dependent-bounds}
For every $i\in\{1,\ldots,M\}$, there exists $E_i \in \mathbb R_+^{m \times n}$ such that
\begin{enumerate}[label=\textit{(\roman*)}]
\item $A_i\geq |B_i|E_i$;
\item $s\geq E_i^\top|r|$;
\item the algebraic equation
\begin{equation}
p_i=s+A_i^\top p_i
-E_i^\top\bigl|r+B_i^\top p_i\bigr|
\label{eq:model-dependent-algebraic-equation}
\end{equation}
admits a solution $p_i\in\mathbb R_+^n$.
\end{enumerate}
\end{assumption}
\begin{remark}
The three items imply the following.
\begin{enumerate}[label=\textit{(\roman*)},leftmargin=*,itemsep=2pt]
\item Every $|u|\leq E_i x$ preserves positivity:
$
A_i x+B_i u
\geq(A_i-|B_i|E_i)x\geq0.
$
\item The corresponding stage cost is nonnegative:
$
s^\top x+r^\top u
\geq(s-E_i^\top|r|)^\top x\geq0.
$
\item The solution $p_i$ is the model-based optimal cost vector, and the
corresponding optimal gain is
$
K_i
\coloneqq
-\operatorname{diag}\!\bigl(
\operatorname{sgn}(r+B_i^\top p_i)
\bigr)E_i
\eqqcolon
\Lambda_iE_i,
$
where $\Lambda_i$ is the corresponding diagonal sign matrix.
They satisfy $A_i+B_iK_i\geq0$, $\rho(A_i+B_iK_i)<1$ and
\[
p_i=s+K_i^\top r+(A_i+B_iK_i)^\top p_i.
\]
Thus, $p_i^\top x$ is the model-based optimal cost, while $K_i x$ is
the corresponding stabilizing optimal input. The disturbance penalty
vector $\gamma$ must satisfy $\gamma\geq p_i$ componentwise for every
model $i$.
\end{enumerate}
\end{remark}
When the true model is known to be model $i$, the controller uses the
model-based gain $K_i$ and input bound $E_i$. If the true model is only known
to belong to $\mathcal C$, however, the input must preserve positivity for
every model in $\mathcal C$. For $q=1,\ldots,m$ and
$\ell=1,\ldots,n$, define
\begin{equation}
	\bigl[E(\mathcal C)\bigr]_{q\ell}
	\coloneqq
	\min_{i:(A_i,B_i)\in\mathcal C}[E_i]_{q\ell},
	\quad
	\mathcal U(x,\mathcal C)
	\coloneqq\left\{u:\ |u|\leq E(\mathcal C)x\right\}.
	\label{eq:consistency-dependent-gain}
\end{equation}
Thus, each entry of $E(\mathcal C)$ is the smallest corresponding entry among
the bounds associated with the models in $\mathcal C$. Consequently, under Assumption~\ref{ass:model-dependent-bounds}, every
input in $\mathcal U(x,\mathcal C)$ preserves positivity for all remaining
models. As models are eliminated, the bound may relax and permit less cautious
control actions.
To exploit this consistency-dependent bound, associate with every
$(A_k,B_k)\in\mathcal C$ a gain $K_k(\mathcal C)$ satisfying
$
|K_k(\mathcal C)|\leq E(\mathcal C).
$
For $x\geq0$, this condition implies
$|K_k(\mathcal C)x|\leq E(\mathcal C)x$, so
$K_k(\mathcal C)x\in\mathcal U(x,\mathcal C)$. When
$\mathcal C=\{(A_i,B_i)\}$, set $K_i(\mathcal C)=K_i$.
The history vector selects the model and its associated admissible gain
according to
\begin{equation}
\begin{aligned}
k
&\in
\arg\min_{\substack{i:(A_i,B_i)\in\mathcal C}}
\gamma^\top[\,I\ {-}A_i\ {-}B_i\,]Z,\\
u
&=\overbar\eta(x,\mathcal C,Z)
=K_k(\mathcal C)x.
\end{aligned}
\label{eq:finite-family-online-controller}
\end{equation}
The positivity guarantee established above can also be verified directly.
For every $(A_i,B_i)\in\mathcal C$,
\begin{equation*}
A_i+B_iK_k(\mathcal C)
\geq A_i-|B_i|E(\mathcal C)
\geq A_i-|B_i|E_i
\geq0.
\label{eq:consistency-dependent-safety}
\end{equation*}
Hence, the selected gain preserves positivity for every model consistent with
the data, even when the selected model $k$ is not the true model. As
$\mathcal C$ shrinks, $E(\mathcal C)$ increases, which may permit less
conservative gain selection. The following theorem provides
performance guarantees for controllers of the
form~\eqref{eq:finite-family-online-controller}.

\begin{theorem}
\label{thm2}
Suppose Assumptions~\ref{ass:model-uncertainty}--
\ref{ass:model-dependent-bounds} hold and that there exist vectors
$p_{ij}\in\mathbb R_+^n$ satisfying
$0\leq p_i\leq p_{ii}\leq\gamma$ for
$i\in\{1,\ldots,M\}$ and
$0\leq p_{ij}=p_{ji}\leq\gamma$ for $i\neq j$.
For every consistency set
$\mathcal C\subseteq\mathcal M$, every
$(A_i,B_i),(A_j,B_j),(A_k,B_k)\in\mathcal C$, except when
$i\neq j=k$, and every $x\in\mathbb R_+^n$, suppose that
\begin{multline}
p_{jk}^\top x
\geq
s^\top x+r^\top K_k(\mathcal C)x
+p_{ij}^\top
\max_{\ell\in\{i,j\}}
\bigl(A_\ell+B_\ell K_k(\mathcal C)\bigr)x
-\tfrac12\gamma^\top
\left|
\bigl[(A_i-A_j)
+(B_i-B_j)K_k(\mathcal C)\bigr]x
\right|.
\label{eq:finite-family-inequality}
\end{multline}
For every $i,j\in\{1,\ldots,M\}$, let
\begin{equation*}
V_{ij}(x,Z)
\coloneqq
p_{ij}^\top x
-\tfrac12\gamma^\top
\Bigl(
[\,I\ {-}A_i\ {-}B_i\,]
+[\,I\ {-}A_j\ {-}B_j\,]
\Bigr)Z,
\label{eq:pairwise-value}
\end{equation*}
and define
\begin{equation}
\overline V(x,\mathcal C,Z)
\coloneqq
\begin{cases}
p_i^\top x
-\gamma^\top[\,I\ {-}A_i\ {-}B_i\,]Z,
&\hspace{-4mm}\mathcal C=\{(A_i,B_i)\},\\[1mm]
\displaystyle
\max_{\substack{i,j:
(A_i,B_i),(A_j,B_j)\in\mathcal C}}
V_{ij}(x,Z),
&|\mathcal C|>1.
\end{cases}
\label{eq:finite-family-value}
\end{equation}
Then the controller in~\eqref{eq:finite-family-online-controller} is
admissible and satisfies the Bellman inequality, namely,
$\mathcal F\overline V
\leq\mathcal F_{\overbar\eta}\overline V
\leq\overline V$.
Moreover, the induced policy $\overbar\mu$, obtained by substituting
$\mathcal C_t$ and $Z_t$ from \eqref{eq:Ct-history} into
\eqref{eq:finite-family-online-controller}, achieves the cost bound
$J_{\overbar\mu}(x_0)
\leq\max_{i,j\in\{1,\ldots,M\}}p_{ij}^\top x_0$.
\end{theorem}
{\em Proof.}
See Appendix~\ref{thm2-proof}. \qedblack
\begin{remark}[Reduced model search]
The controller minimizes
$\gamma^\top[\,I\ {-}A_i\ {-}B_i\,]Z_t$
only over the models in the current consistency set. Since
$\mathcal C_{t+1}\subseteq\mathcal C_t$, models incompatible with an observed
transition are permanently removed. The same consistency set determines
$E(\mathcal C_t)$, which can increase componentwise as uncertainty decreases,
thereby enlarging the admissible input set.
\end{remark}

\begin{remark}[Common input-bound case] \label{Remark5}
Suppose that all model-dependent bounds coincide,
$E_1=\cdots=E_M=E$. Then $E(\mathcal C)=E$ for every nonempty
$\mathcal C\subseteq\mathcal M$. If, in addition,
$K_k(\mathcal C)=K_k$ whenever $(A_k,B_k)\in\mathcal C$, then
\eqref{eq:finite-family-inequality} becomes
\begin{align*}
p_{jk}^\top x
\geq
s^\top x+r^\top K_kx
+p_{ij}^\top
\max_{\ell\in\{i,j\}}
(A_\ell+B_\ell K_k)x
-\tfrac12\gamma^\top
\left|
\bigl[(A_i-A_j)+(B_i-B_j)K_k\bigr]x
\right|.
\label{eq:common-bound-condition}
\end{align*}
This condition coincides with that of
\cite[Theorem~3]{Bencherki2026}. Hence, the present theorem recovers the
common bound result as a special case using the same value function ansatz.
\end{remark}

We next specialize \eqref{eq:finite-family-inequality} to two models to
clarify the meaning of the inequalities in the sufficient conditions.

\begin{example}[Two-model set specialization]
Let $\mathcal M=\{(A_1,B_1),(A_2,B_2)\}$ and define
\begin{equation*}
E_{12}\coloneqq\min\{E_1,E_2\},
\qquad
|K_k^{12}|\leq E_{12},
\quad k\in\{1,2\},
\label{eq:two-model-safe-gains}
\end{equation*}
where the minimum is elementwise and
$K_k^{12}\coloneqq K_k(\mathcal M)$. These gains preserve positivity for
both models because
\begin{equation*}
A_i+B_iK_k^{12}
\geq A_i-|B_i|E_{12}
\geq A_i-|B_i|E_i
\geq0,
\quad i,k\in\{1,2\}.
\label{eq:two-model-positivity}
\end{equation*}
Let $0\leq p_{11},p_{12}=p_{21},p_{22}\leq\gamma$. The sufficient conditions in~\eqref{eq:finite-family-inequality}
have the following interpretations.
Let $k$ denote the model selected by the controller and $\mathcal C_+$ the
consistency set after the next observation. The conditions in~\eqref{eq:finite-family-inequality} correspond to
three cases:
\begin{enumerate}[label=\textit{(\roman*)}]

\item The model is already known, i.e.,
$\mathcal C=\{(A_i,B_i)\}$
\begin{equation}
p_{i}\geq
s+K_i^\top r+(A_i+B_iK_i)^\top p_{i},
\qquad i\in\{1,2\}.
\label{eq:two-model-singleton-all}
\end{equation}

\item The next observation isolates model $i$, i.e.,
$\mathcal C_+=\{(A_i,B_i)\}$
\begin{equation*}
p_{ik}\geq
s+(K_k^{12})^\top r
+(A_i+B_iK_k^{12})^\top p_{ii},
\qquad i,k\in\{1,2\}.
\label{eq:two-model-diagonal-all}
\end{equation*}

\item Both models remain consistent, i.e.,
$\mathcal C_+=\mathcal C$
\begin{multline}
p_{12}^\top x
\geq
s^\top x+r^\top K_k^{12}x
+p_{12}^\top
\max_{\ell\in\{1,2\}}
(A_\ell+B_\ell K_k^{12})x
-\tfrac12\gamma^\top
\left|
\bigl[(A_1-A_2)+(B_1-B_2)K_k^{12}\bigr]x
\right|,
\label{eq:two-model-ambiguity}
\end{multline}
for every $k\in\{1,2\}$ and $x\in\mathbb R_+^n$.
\end{enumerate}
Item~\textit{(i)} applies after the true model has been identified,
item~\textit{(ii)} covers transitions from model uncertainty to identification
of the true model, and item~\textit{(iii)} addresses the case in which both
candidate models remain consistent with the data. Before the true model is
identified, the controller applies $K_k^{12}$ associated with the candidate
model $k$ selected from the observed history. The bound
$|K_k^{12}|\leq E_{12}$ ensures that this gain preserves positivity for both
candidate models, even when the selected model is not the true one. Stability
is required only for the selected model, as enforced by the case $i=k$ in
item~\textit{(ii)}; the effect of applying $K_k^{12}$ when the other candidate
is the true model is captured by the cross cost vector $p_{12}$. These
selected-model stability conditions are additional to those
in~\cite{Bencherki2026}. When $E_1=E_2$ and $K_k^{12}=K_k$, they coincide with
item~\textit{(i)}, thereby recovering the common bound result stated in Remark~\ref{Remark5}.
\end{example}
The results in this section apply to finite model sets.
A central contribution of this paper is to identify
structured continuum model classes for which the required
information admits a fixed, finite-dimensional recursive
representation, even when $\mathcal M$ is uncountable.
These classes are developed in the next section.
\section{History compression for continuum model classes}
\label{sec:history-compression}
The formulation in \cite{Bencherki2026} considers positivity-preserving
disturbances satisfying $w_t\geq-(Ax_t+Bu_t)$, where
$Ax_t+Bu_t\geq0$ for $u_t\in\mathcal U(x_t,\mathcal{C}_t)$. With $v_t=x_{t+1}$, this
condition is equivalent to $v_t\geq0$, independently of the true pair $(A_\star,B_\star) \in \mathcal{M}$ selected by the adversary.
Consequently, every candidate model $(A,B)$ can explain any nonnegative observation
through $w_t^{(A,B)}=v_t-Ax_t-Bu_t$ and therefore remains feasible.
Because
$
w_t^{(A,B)}=v_t-Ax_t-Bu_t
$
may have either sign, its cumulative absolute-value penalty cannot generally
be computed from a common history vector. The formulation in
\cite{Bencherki2026} therefore maintains a separate historical variable for
each $(A,B)\in\mathcal M$.

In the present work, the controller selects
$u_t\in\mathcal U(x_t,\mathcal C_t)$, which preserves positivity for every
model in $\mathcal C_t$. For a candidate $(A,B)$, the condition $w_t\geq0$
is equivalent to $v_t\geq Ax_t+Bu_t$. An observation may therefore
eliminate incompatible models, yielding
$\mathcal C_{t+1}=\{(A,B)\in\mathcal C_t:v_t\geq Ax_t+Bu_t\} \subseteq \mathcal{C}_t$. For every
$(A,B)\in\mathcal C_t$, the disturbances reconstructed from the observed
history are nonnegative.
Hence,
\begin{equation}
\sum_{\tau=0}^{t-1}\gamma^\top w_\tau^{(A,B)}
=
\gamma^\top[\,I\ {-}A\ {-}B\,]
\overbrace{
\sum_{\tau=0}^{t-1}
\begin{bmatrix}
v_\tau^\top & x_\tau^\top & u_\tau^\top
\end{bmatrix}^{\!\top}
}^{Z_t}.
\label{eq:continuum-common-history}
\end{equation}
Thus, the same fixed-dimensional vector $Z_t$ recovers the accumulated
disturbance penalty for every model in $\mathcal C_t$. The two disturbance
classes are compared in Table~\ref{tab:compression-comparison}.

\begin{table*}[!t]
	\centering
	\vspace{-2mm}
	\caption{Comparison of learning and history compression under the two
		disturbance classes.}
	\label{tab:compression-comparison}
	\small
	\renewcommand{\arraystretch}{1.25}
	\setlength{\tabcolsep}{7pt}
	\begin{tabularx}{\textwidth}{
			@{}>{\raggedright\arraybackslash}p{0.14\textwidth}
			>{\raggedright\arraybackslash}X
			>{\raggedright\arraybackslash}X@{}}
		\toprule
		\textbf{Property}
		&
		\centering\arraybackslash
		\textbf{Positivity-preserving disturbances in~\cite{Bencherki2026,Bencherki2026b}}
		&
		\centering\arraybackslash
		\textbf{Nonnegative disturbances in the current work}
		\\
		\midrule
		\textbf{Learning}
		&
		No feasibility-based learning: every $v_t\geq0$ is feasible for every
		model $(A,B)\in\mathcal M$.
		&
		Feasibility-based learning: observations eliminate models violating
		$v_t\geq Ax_t+Bu_t$, shrinking $\mathcal C_t\subseteq\mathcal M$.
		\\
		\textbf{Information state}
		&
		$x_t$ together with a model-dependent history variable $z_t^{(A,B)}$
		for each $(A,B)\in\mathcal M$.
		&
		The compressed information state $(x_t,\mathcal C_t,Z_t)$, with one
		common history vector $Z_t$.
		\\
		\textbf{Memory}
		&
The representation maintains one model-specific historical variable for each candidate model.
		&
		$Z_t$ has fixed dimension; structured consistency sets may also admit
		fixed-dimensional representations.
		\\
		\textbf{Input bound}
		&
		A common bound $E$ preserves positivity for every model in $\mathcal M$
		and remains fixed throughout operation; consequently, the search for
		stabilizing gains is restricted to gains satisfying $|K|\leq E$.
		&
		The consistency-dependent bound $E(\mathcal C_t)$ preserves positivity
		for every remaining model and allows the gain magnitude to vary with
		$\mathcal C_t$. As the consistency set contracts, the bound may increase,
		permitting less cautious control actions.
		\\
		\bottomrule
	\end{tabularx}
\end{table*}

Although $Z_t$ has fixed dimension, $\mathcal C_t$ need not. For structured
continuum model classes, further compression is possible. Suppose that each
$(A,B)\in\mathcal M$ admits a bound $E_{(A,B)}$ satisfying
$A\geq|B|E_{(A,B)}$ and an associated model-based gain
$K_{(A,B)}=\Lambda_{(A,B)}E_{(A,B)}$. For every nonempty
$\mathcal C\subseteq\mathcal M$, define
\begin{equation*}
E(\mathcal C)
\coloneqq
\inf_{(A,B)\in\mathcal C}E_{(A,B)},
\qquad
\mathcal U(x,\mathcal C)
\coloneqq
\{u:|u|\leq E(\mathcal C)x\},
\label{eq:continuum-safe-bound}
\end{equation*}
where the infimum is taken entrywise. The finite-model CE law from Section~\ref{sec:finite-model-certificate} extends to continuum model classes as
\begin{equation}
	u=\overbar\eta(x,\mathcal C,Z)
	\coloneqq\Lambda_{(\widehat A,\widehat B)}E(\mathcal C)x,
	\quad
	(\widehat A,\widehat B)\in
	\operatorname*{arg\,min}_{(A,B)\in\mathcal C}
	\gamma^\top[\,I\ {-}A\ {-}B\,]Z.
	\label{eq:continuum-CE-law}
\end{equation}
The history selects a candidate model and the sign structure of its
model-based feedback gain, while $\mathcal C$ determines a safe input bound that preserves positivity for every remaining model. As $\mathcal C$
contracts, $E(\mathcal C)$ can increase, allowing progressively less cautious
inputs. We next present a continuum model class for which both
$\mathcal C_t$ and $E(\mathcal C_t)$ admit finite-dimensional
representations.

\begin{example}[Interval uncertainty in the input matrix]
	\label{example1}
	Consider
	\begin{equation}
		x_{t+1}=Ax_t+B_\star u_t+w_t,
		\qquad
		x_0\in\mathbb R_+^n,\quad w_t\geq0,
	\end{equation}
	where $A\in\mathbb R_+^{n\times n}$ is known, $u_t\in\mathbb R$, and
	\begin{equation}
		(A,B_\star)\in\mathcal C_0=\mathcal M
		\coloneqq
		\left\{(A,B):B\in[\underline B_0,\overline B_0]\right\},
		\qquad
		(A,\mathbf0)\notin\mathcal C_0.
	\end{equation}
	Here $[\underline B_0,\overline B_0]\subset\mathbb R^{n\times1}$.
	The interval and all minimum and maximum operations involving its endpoints
	are understood elementwise.
	Suppose that
	$\mathcal C_t=\{(A,B):B\in[\underline B_t,\overline B_t]\}$ and define
	\begin{equation}
		\begin{aligned}
			\beta_t
			&\coloneqq
			\max\left\{
			|\underline B_t|,
			|\overline B_t|
			\right\},\\
			[E(\mathcal C_t)]_\ell
			&\coloneqq
			\min_{\substack{q\in\{1,\ldots,n\}\\
					[\beta_t]_q>0}}
			\tfrac{[A]_{q\ell}}{[\beta_t]_q},
			\qquad
			\ell\in\{1,\ldots,n\}.
		\end{aligned}
		\label{eq:interval-adaptive-bound}
	\end{equation}
	It holds for every $(A,B)\in\mathcal C_t$ that
	\[
	|[B]_q|[E(\mathcal C_t)]_\ell
	\leq
	[\beta_t]_q[E(\mathcal C_t)]_\ell
	\leq
	[A]_{q\ell},
	\]
	and hence
	$
	A\geq|B|E(\mathcal C_t).
	$
	Therefore, the input set
	\begin{equation*}
		\mathcal U(x_t,\mathcal C_t)
		\coloneqq
		\left\{
		u\in\mathbb R:
		|u|\leq E(\mathcal C_t)x_t
		\right\}
		\label{eq:interval-adaptive-input-set}
	\end{equation*}
	preserves positivity for every $(A,B)\in\mathcal C_t$, since
	$
	Ax_t+Bu_t
	\geq
	\bigl(A-|B|E(\mathcal C_t)\bigr)x_t
	\geq0.
	$
	The condition $(A,\mathbf0)\notin\mathcal C_0$ excludes a model with no input
	influence.
	For $u_t\neq0$, let
	$
	Q_t\coloneqq\tfrac{v_t-Ax_t}{u_t}.
	$
	The interval endpoints are updated elementwise according to
	\begin{equation}
		(\underline B_{t+1},\overline B_{t+1})
		=
		\begin{cases}
			\bigl(
			\underline B_t,
			\min\{\overline B_t,Q_t\}
			\bigr),
			&u_t>0,\\[1mm]
			\bigl(
			\max\{\underline B_t,Q_t\},
			\overline B_t
			\bigr),
			&u_t<0.
		\end{cases}
		\label{eq:interval-consistency-update}
	\end{equation}
	Thus,
	$
	\mathcal C_{t+1}
	=
	\{(A,B):B\in[\underline B_{t+1},\overline B_{t+1}]\}.
	$
	If $u_t=0$, feasibility requires $v_t\geq Ax_t$, in which case
	$\mathcal C_{t+1}=\mathcal C_t$. Moreover,
	\[
	v_t
	=
	Ax_t+B_\star u_t+w_t
	\geq
	Ax_t+B_\star u_t,
	\]
	so $(A,B_\star)\in\mathcal C_{t+1}$. Hence, the consistency set remains
	nonempty and is represented exactly by its two endpoints.
	The CE estimate is
	\begin{equation}
		\widehat B_t
		\in
		\arg\min_{B:(A,B)\in\mathcal C_t}
		\gamma^\top[\,I\ {-}A\ {-}B\,]Z_t.
		\label{eq:interval-CE-estimate}
	\end{equation}
Let
$\Lambda_{(A,\widehat B_t)}\in\{-1,0,1\}$ denote the model-based
feedback direction. The adaptive controller is
\begin{equation}
	u_t=K_{(A,\widehat B_t)}x_t,\qquad
	K_{(A,\widehat B_t)}
	\coloneqq
	\Lambda_{(A,\widehat B_t)}E(\mathcal C_t).
	\label{eq:interval-adaptive-controller}
\end{equation}
Since
$\mathcal C_{t+1}\subseteq\mathcal C_t$, we have
$
\beta_{t+1}\leq\beta_t
\Longrightarrow
E(\mathcal C_{t+1})
\geq
E(\mathcal C_t)
$
elementwise. The admissible input range may therefore expand as
inconsistent models are eliminated. Establishing stabilization
guarantees for the CE controller under continuum model sets,
such as the one considered here, is deferred to future work.
\end{example}

\begin{example}[Interval vector uncertainty in the state matrix]
	Consider the system
	\begin{align*}
		x_{t+1}=(A_0+a_\star c^\top)x_t+Bu_t+w_t,
	\end{align*}
	where \(A_0\in\mathbb R^{n\times n}\),
	\(B\in\mathbb R^{n\times m}\), and \(c\in\mathbb R^n\)
	are known, with no sign restriction on \(c\).
	The unknown vector \(a_\star\in\mathbb R^n\) satisfies
	\(\underline a_0\leq a_\star\leq\overline a_0\).
	The states and disturbances satisfy \(x_t,w_t\in\mathbb R_+^n\).
	Thus, with \(A_\star=A_0+a_\star c^\top\),
	\[
	(A_\star,B)\in\mathcal C_0=\mathcal M
	=\{(A_0+ac^\top,B):
	\underline a_0\leq a\leq\overline a_0\}.
	\]
	Assume \(|u_t|\leq Ex_t\), where
	\(E\in\mathbb R_+^{m\times n}\), and
	\(A_0+ac^\top\geq|B|E\) for every
	\(a\in[\underline a_0,\overline a_0]\).
	It follows that \((A_0+ac^\top)x_t+Bu_t\geq0\) for every
	admissible input and candidate matrix.
	
	Suppose
	\(\mathcal C_t=\{(A_0+ac^\top,B):
	\underline a_t\leq a\leq\overline a_t\}\).
	A candidate \(a\) is consistent with the observation
	\(v_t=x_{t+1}\) if
	\[
	v_t\geq(A_0+ac^\top)x_t+Bu_t
	=A_0x_t+a(c^\top x_t)+Bu_t.
	\]
	Hence \(a(c^\top x_t)\leq v_t-A_0x_t-Bu_t\).
	If \(c^\top x_t=0\), this condition is independent of \(a\),
	so a feasible observation leaves \(\mathcal C_t\) unchanged.
	Otherwise, define
	\(q_t\coloneqq(v_t-A_0x_t-Bu_t)/(c^\top x_t)\).
	Consistency then requires \(a\leq q_t\) when \(c^\top x_t>0\),
	and \(a\geq q_t\) when \(c^\top x_t<0\).
	The interval endpoints therefore satisfy
	\[
	(\underline a_{t+1},\overline a_{t+1})
	=
	\begin{cases}
		\bigl(\underline a_t,\min\{\overline a_t,q_t\}\bigr),
		& c^\top x_t>0,\\[1mm]
		\bigl(\max\{\underline a_t,q_t\},\overline a_t\bigr),
		& c^\top x_t<0,\\[1mm]
		\bigl(\underline a_t,\overline a_t\bigr),
		& c^\top x_t=0,
	\end{cases}
	\]
	where minimum and maximum are componentwise, and
	\[
	\mathcal C_{t+1}
	=\{(A_0+ac^\top,B):
	\underline a_{t+1}\leq a\leq\overline a_{t+1}\}.
	\]
	The observation generated by \(a_\star\) satisfies
	\(a_\star\leq q_t\) when \(c^\top x_t>0\) and
	\(a_\star\geq q_t\) when \(c^\top x_t<0\).
	Hence the true parameter remains in the updated interval.
	Thus, the consistency set can be maintained by updating only
	the two endpoint vectors \(\underline a_t\) and \(\overline a_t\),
	which exactly represent the entire set of consistent models.
	The control input \(u_t\in\mathbb R^m\) may be multidimensional.
	\label{example3}
\end{example}

\section{Numerical Experiment}
\label{sec:numerical-experiments}
We illustrate the proposed approach using a two-state positive system with
a scalar input and interval uncertainty in the input
matrix.\footnote{Implementation details and code are available in the
	\href{https://github.com/Fethi-Bencherki/minimax-adaptive-positive-control-nonnegative-disturbances}{{\tt GitHub repository}}.}
Consider
\begin{equation*}
	x_{t+1}=Ax_t+B_\star u_t+w_t,
	\quad
	A=\begin{bmatrix}0.45&0.20\\0.20&0.50\end{bmatrix},
	\
	x_0=\begin{bmatrix}1\\0.7\end{bmatrix}.
	\label{eq:numerical-system}
\end{equation*}
The matrix $A$ is known and the fixed model pair belongs to
\begin{equation}
	(A,B_\star)\in\mathcal C_0=\mathcal M
	=\left\{(A,B):B\in[\underline B_0,\overline B_0]\right\},
	\quad
	[\underline B_0,\overline B_0]
	=
	\left[
	\begin{bmatrix}
		0.05\\
		-0.80
	\end{bmatrix},
	\begin{bmatrix}
		0.45\\
		0.10
	\end{bmatrix}
	\right].
	\label{eq:numerical-initial-set}
\end{equation}
The true input matrix is
$
B_\star=\begin{bmatrix}0.20\\-0.40\end{bmatrix}
$
and remains fixed throughout the experiment. Notice that
$(A,\mathbf0)\notin\mathcal C_0$.
Applying \eqref{eq:interval-adaptive-bound} to $\mathcal C_0$ and to the
singleton $\{(A,B_\star)\}$ gives
\begin{equation}
	E(\mathcal C_0)
	=
	\begin{bmatrix}0.25&0.44\end{bmatrix},
	\;\;
	E(\{(A,B_\star)\})
	=
	\begin{bmatrix}0.50&1.00\end{bmatrix},
	\label{eq:numerical-bounds}
\end{equation}
where
$
[E(\{(A,B_\star)\})]_\ell
=
\min_{q:\,|[B_\star]_q|>0}
[A]_{q\ell}/|[B_\star]_q|
$.
Thus, $E(\{(A,B_\star)\})$ is the full model-specific bound available when
$B_\star$ is known.
We take $s=\mathbf1$, $r=0$, and $\gamma=6\cdot\mathbf1$. Since the input is
scalar and $\gamma>0$, the objective in
\eqref{eq:interval-CE-estimate} is affine in $B$. Its minimum over
$\mathcal C_t$ is therefore attained at an endpoint:
\begin{equation*}
	\widehat B_t=\overline B_t
	\ \text{if }\sum_{\tau<t}u_\tau\geq0,
	\qquad
	\widehat B_t=\underline B_t\ \text{otherwise}.
	\label{eq:numerical-endpoint-rule}
\end{equation*}
Thus, the CE estimate is computed directly from the interval endpoints,
without discretizing the continuum $\mathcal C_t$. We apply the adaptive
controller in \eqref{eq:interval-adaptive-controller}. The baseline knows
$B_\star$ and uses the full model-specific bound $E(\{(A,B_\star)\})$.
The disturbances are sampled independently according to
$[w_t]_1,[w_t]_2\sim\operatorname{Unif}[0,0.2].$
The same disturbance realization is applied to the adaptive controller and
to the baseline controller that knows $B_\star$.

\begin{figure}[H]
\centering
\begin{subfigure}[t]{0.47\textwidth}
\centering
\begin{tikzpicture}

\definecolor{chocolate213940}{RGB}{213,94,0}
\definecolor{darkcyan0114178}{RGB}{0,114,178}
\definecolor{darkgrey176}{RGB}{176,176,176}

\begin{axis}[
    width=0.85\linewidth,
height=0.45\linewidth,
legend cell align={left},
legend columns=2,
legend style={
  fill opacity=0.8,
  draw opacity=1,
  text opacity=1,
        at={(0.5,1.5)},
  anchor=north,
  draw=none
},
tick align=outside,
tick pos=left,
width=\linewidth,
x grid style={darkgrey176},
xlabel={$t$},
xmajorgrids,
xmin=-2.45, xmax=51.45,
xminorgrids,
xtick style={color=black},
y grid style={darkgrey176},
ylabel={$x_{1,t}$},
ymajorgrids,
ymin=0.0804712266117758, ymax=1.04378708444706,
yminorgrids,
ytick style={color=black}
]
\addplot [semithick, darkcyan0114178, mark=*, mark size=1.5, mark options={solid}]
table {%
0 1
1 0.618341579254267
2 0.659476809102094
3 0.334431766264611
4 0.573987528149887
5 0.331910391211762
6 0.423580059500331
7 0.187579152132066
8 0.29180242171136
9 0.323829824098422
10 0.376686646674553
11 0.33024690041924
12 0.36517721095883
13 0.261064583798059
14 0.292078191305248
15 0.311548395477916
16 0.232785914688317
17 0.357195809710755
18 0.154870219548342
19 0.285285462280173
20 0.254378085874586
21 0.450803184438618
22 0.253676551852544
23 0.33383830302996
24 0.21959340066431
25 0.138629556639119
26 0.37352835901009
27 0.242906917725712
28 0.316117494437013
29 0.369472002566257
30 0.269601699842818
31 0.171592149281127
32 0.13434156210233
33 0.126472073478769
34 0.232833734721191
35 0.334461923285367
36 0.237170854986971
37 0.390474170881773
38 0.363522124383562
39 0.396223957523179
40 0.4205201397024
41 0.299989042583177
42 0.266019946108399
43 0.385470894377327
44 0.285639346397518
45 0.341074700872212
46 0.418542553320307
47 0.301506882332463
48 0.280848882559093
49 0.355380072838225
};
\addlegendentry{Adaptive}
\addplot [semithick, chocolate213940, dashed, mark=square*, mark size=1.5, mark options={solid}]
table {%
0 1
1 0.970563801476489
2 0.693424985716908
3 0.493231728689656
4 0.491110372594023
5 0.420451762046888
6 0.357776397619738
7 0.284170889949681
8 0.219894902381147
9 0.392718244542741
10 0.330196859936258
11 0.413681426372212
12 0.307058503162692
13 0.350790950017216
14 0.238075319920281
15 0.391755425593408
16 0.395797140955422
17 0.271260428041659
18 0.21981204086789
19 0.215754842874305
20 0.300397710573171
21 0.400230618114467
22 0.375954196391547
23 0.279092920456121
24 0.286941441586351
25 0.229043553776555
26 0.307337773773084
27 0.184413689257547
28 0.277609186622797
29 0.348718663965358
30 0.258596208151951
31 0.165148393709373
32 0.130128726175931
33 0.124258311058834
34 0.232314782627836
35 0.335094531843448
36 0.238367313917491
37 0.391308392940412
38 0.365667208902948
39 0.397083077729145
40 0.421573870825187
41 0.301666758406652
42 0.266446221569637
43 0.385959130630075
44 0.287407443447369
45 0.342145621718233
46 0.41937212788493
47 0.303269433491097
48 0.281269532147348
49 0.355302333814925
};
\addlegendentry{Known model}
\end{axis}

\end{tikzpicture}
\caption{First state.}
\label{fig:state-x1}
\end{subfigure}\hfill
\begin{subfigure}[t]{0.47\textwidth}
\centering
\begin{tikzpicture}

\definecolor{chocolate213940}{RGB}{213,94,0}
\definecolor{darkcyan0114178}{RGB}{0,114,178}
\definecolor{darkgrey176}{RGB}{176,176,176}

\begin{axis}[
    width=0.85\linewidth,
height=0.45\linewidth,
legend cell align={left},
legend columns=2,
legend style={
  fill opacity=0.8,
  draw opacity=1,
  text opacity=1,
       at={(0.5,1.5)},
  anchor=north,
  draw=none
},
tick align=outside,
tick pos=left,
width=\linewidth,
x grid style={darkgrey176},
xlabel={$t$},
xmajorgrids,
xmin=-2.45, xmax=51.45,
xminorgrids,
xtick style={color=black},
y grid style={darkgrey176},
ylabel={$x_{2,t}$},
ymajorgrids,
ymin=-0.0314994400210299, ymax=0.943526871184688,
yminorgrids,
ytick style={color=black}
]
\addplot [semithick, darkcyan0114178, mark=*, mark size=1.5, mark options={solid}]
table {%
0 0.7
1 0.89920749340261
2 0.494197483588695
3 0.616971561377964
4 0.29652160039937
5 0.588483866915439
6 0.348306853023618
7 0.434913438868619
8 0.295447097471261
9 0.513100773595268
10 0.239571697529563
11 0.39126568583255
12 0.243863402188237
13 0.380256072966709
14 0.236603901379809
15 0.486467858452998
16 0.5484669299404
17 0.230581561988814
18 0.460473526656336
19 0.161075705016679
20 0.367302963213166
21 0.241233103718125
22 0.410555210720214
23 0.174092196937876
24 0.320390593577156
25 0.483807363036332
26 0.109193659526711
27 0.101965117337168
28 0.113820988228721
29 0.0311509690587047
30 0.0547857780872706
31 0.0564025942001469
32 0.109308119081245
33 0.201562842295685
34 0.135981570639096
35 0.0530278336698515
36 0.16079598372989
37 0.14323294085264
38 0.0308553876396747
39 0.0210799232345342
40 0.0894720166879624
41 0.102854231166974
42 0.147942489029194
43 0.147746252144417
44 0.218916904457386
45 0.105225841502537
46 0.131189509263181
47 0.123457598985481
48 0.0552298082479635
49 0.176812264260303
};
\addlegendentry{Adaptive}
\addplot [semithick, chocolate213940, dashed, mark=square*, mark size=1.5, mark options={solid}]
table {%
0 0.7
1 0.194763048958165
2 0.179745574803512
3 0.053707912897506
4 0.130018586546235
5 0.170383541665913
6 0.20101939197887
7 0.0367872202244625
8 0.187199450951323
9 0.198802779575451
10 0.12546033887707
11 0.0705576457272365
12 0.163241560234
13 0.0643131041938297
14 0.158959975095225
15 0.196771106042424
16 0.0499451333637156
17 0.133095057301805
18 0.148323109962927
19 0.0906375722558506
20 0.135385712984793
21 0.184274298329335
22 0.059107173835731
23 0.101785137882422
24 0.0603982377519423
25 0.143069093533012
26 0.0343657842989141
27 0.0787968428077758
28 0.105643605465718
29 0.0229388631781629
30 0.0465771937155271
31 0.0497963385491342
32 0.104700029719201
33 0.197208967190677
34 0.130665848569371
35 0.0464076193882478
36 0.153140742429293
37 0.135990884388399
38 0.0209650120475868
39 0.0128199377610482
40 0.0808755997428146
41 0.092921128504456
42 0.139995633116274
43 0.140086512324655
44 0.20902400439792
45 0.0961252965118122
46 0.122517882572906
47 0.113040564668051
48 0.0469519844602647
49 0.169980415445054
};
\addlegendentry{Known model}
\end{axis}

\end{tikzpicture}
\caption{Second state.}
\label{fig:state-x2}
\end{subfigure}
\par\vspace{7mm}
\begin{subfigure}[t]{0.47\textwidth}
\centering
\begin{tikzpicture}

\definecolor{chocolate213940}{RGB}{213,94,0}
\definecolor{darkcyan0114178}{RGB}{0,114,178}
\definecolor{darkgrey176}{RGB}{176,176,176}

\begin{axis}[
    width=0.85\linewidth,
height=0.45\linewidth,
legend cell align={left},
legend columns=2,
legend style={
  fill opacity=0.8,
  draw opacity=1,
  text opacity=1,
     at={(0.5,1.5)},
  anchor=north,
  draw=none
},
tick align=outside,
tick pos=left,
width=\linewidth,
x grid style={darkgrey176},
xlabel={$t$},
xmajorgrids,
xmin=-2.45, xmax=51.45,
xminorgrids,
xtick style={color=black},
y grid style={darkgrey176},
ylabel={$u_t$},
ymajorgrids,
ymin=-0.649166666666667, ymax=1.28805555555556,
yminorgrids,
ytick style={color=black}
]
\addplot [semithick, darkcyan0114178, mark=*, mark size=1.5, mark options={solid}]
table {%
0 -0.561111111111111
1 0.598359622177898
2 -0.505610255865107
3 0.508745821997546
4 -0.420109695186217
5 0.49074549122183
6 -0.39839679857542
7 0.357615453295845
8 -0.317334765092615
9 0.468311787518201
10 -0.345327239733321
11 0.435011541182905
12 -0.343249953775986
13 0.395664291347499
14 -0.322187902422856
15 -0.531638247147112
16 0.545774074711208
17 -0.346824948717578
18 0.439850908701788
19 -0.258465156223957
20 0.410074305235576
21 -0.397745694926577
22 0.444487848228972
23 -0.290892776346316
24 -0.356660290892686
25 0.451237020662229
26 0.256744121515226
27 0.208767187924864
28 0.253393816436623
29 0.19741853588261
30 0.175123134974083
31 0.132330009348354
32 0.166746234710402
33 0.252598042637008
34 0.237176832813114
35 0.202775745631445
36 0.263190080387921
37 0.315554601429368
38 0.194439080045762
39 0.199765856001626
40 0.277048434119155
41 0.235464888341892
42 0.263789826512219
43 0.317664384139356
44 0.341458440194198
45 0.256359261460637
46 0.31658610680233
47 0.256120739261823
48 0.180644083436311
49 0.332215850008205
};
\addlegendentry{Adaptive}
\addplot [semithick, chocolate213940, dashed, mark=square*, mark size=1.5, mark options={solid}]
table {%
0 1.2
1 0.68004494969641
2 0.526458067661966
3 0.300323777242334
4 0.375573772843246
5 0.380609422689357
6 0.379907590788739
7 0.178872665199303
8 0.297146902141896
9 0.395161901846822
10 0.290558768845199
11 0.277398358913343
12 0.316770811815346
13 0.239708579202438
14 0.277997635055365
15 0.392648818839128
16 0.247843703841427
17 0.268725271322635
18 0.258229130396872
19 0.198514993693003
20 0.285584568271379
21 0.384389607386569
22 0.247084272031505
23 0.241331598110482
24 0.203868958545118
25 0.25759087042129
26 0.188034671185456
27 0.171003687436549
28 0.244448198777116
29 0.197298195160842
30 0.175875297791503
31 0.13237053540382
32 0.169764392807167
33 0.259338122720094
34 0.246823239883289
35 0.213954885309972
36 0.272324399388039
37 0.331645080858605
38 0.203798616499061
39 0.211361476625621
40 0.291662535155408
41 0.243754507707782
42 0.273218743901093
43 0.333066077639692
44 0.352727726121605
45 0.267198107370929
46 0.332203946515371
47 0.2646752814136
48 0.187586750533939
49 0.347631582352516
};
\addlegendentry{Known model}
\end{axis}

\end{tikzpicture}
\caption{Control input.}
\label{fig:control-input}
\end{subfigure}\hfill
\begin{subfigure}[t]{0.47\textwidth}
\centering
\begin{tikzpicture}

\definecolor{chocolate213940}{RGB}{213,94,0}
\definecolor{darkcyan0114178}{RGB}{0,114,178}
\definecolor{darkgrey176}{RGB}{176,176,176}

\begin{axis}[
    width=0.85\linewidth,
height=0.45\linewidth,
legend cell align={left},
legend columns=2,
legend style={
  fill opacity=0.8,
  draw opacity=1,
  text opacity=1,
        at={(0.5,1.7)},
  anchor=north,
  draw=none
},
tick align=outside,
tick pos=left,
width=\linewidth,
x grid style={darkgrey176},
xlabel={$t$},
xmajorgrids,
xmin=-2.45, xmax=51.45,
xminorgrids,
xtick style={color=black},
y grid style={darkgrey176},
ylabel={$E(\mathcal{C}_t)$},
ymajorgrids,
ymin=0.214039813983556, ymax=1.00516390634533,
yminorgrids,
ytick style={color=black}
]
\addplot [semithick, darkcyan0114178, const plot mark left, mark=*, mark size=1.5, mark options={solid}]
table {%
0 0.25
1 0.321362583399201
2 0.321362583399201
3 0.424922814637528
4 0.424922814637528
5 0.424922814637528
6 0.424922814637528
7 0.452610480949973
8 0.452610480949973
9 0.452610480949973
10 0.452610480949973
11 0.452610480949973
12 0.452610480949973
13 0.452610480949973
14 0.452610480949973
15 0.452610480949973
16 0.452610480949973
17 0.452610480949973
18 0.452610480949973
19 0.452610480949973
20 0.452610480949973
21 0.452610480949973
22 0.452610480949973
23 0.452610480949973
24 0.452610480949973
25 0.452610480949973
26 0.452610480949973
27 0.452610480949973
28 0.452610480949973
29 0.452610480949973
30 0.452610480949973
31 0.452610480949973
32 0.452610480949973
33 0.452610480949973
34 0.452610480949973
35 0.452610480949973
36 0.452610480949973
37 0.452610480949973
38 0.452610480949973
39 0.452610480949973
40 0.452610480949973
41 0.452610480949973
42 0.452610480949973
43 0.452610480949973
44 0.452610480949973
45 0.452610480949973
46 0.452610480949973
47 0.452610480949973
48 0.452610480949973
49 0.452610480949973
};
\addlegendentry{$[E(\mathcal{C}_t)]_1$}
\addplot [semithick, chocolate213940, const plot mark left, mark=square*, mark size=1.5, mark options={solid}]
table {%
0 0.444444444444444
1 0.444444444444444
2 0.594254512725556
3 0.594254512725556
4 0.594254512725556
5 0.594254512725556
6 0.627056187876378
7 0.627056187876378
8 0.627056187876378
9 0.627056187876378
10 0.729781177064725
11 0.729781177064725
12 0.729781177064725
13 0.729781177064725
14 0.802988668712638
15 0.802988668712638
16 0.802988668712638
17 0.802988668712638
18 0.802988668712638
19 0.802988668712638
20 0.802988668712638
21 0.802988668712638
22 0.802988668712638
23 0.802988668712638
24 0.802988668712638
25 0.802988668712638
26 0.802988668712638
27 0.969203720328885
28 0.969203720328885
29 0.969203720328885
30 0.969203720328885
31 0.969203720328885
32 0.969203720328885
33 0.969203720328885
34 0.969203720328885
35 0.969203720328885
36 0.969203720328885
37 0.969203720328885
38 0.969203720328885
39 0.969203720328885
40 0.969203720328885
41 0.969203720328885
42 0.969203720328885
43 0.969203720328885
44 0.969203720328885
45 0.969203720328885
46 0.969203720328885
47 0.969203720328885
48 0.969203720328885
49 0.969203720328885
};
\addlegendentry{$[E(\mathcal{C}_t)]_2$}
\addplot [semithick, darkcyan0114178, dashed]
table {%
-2.45 0.5
51.45 0.5
};
\addlegendentry{$[E(\{(A,B_\star)\})]_1$}
\addplot [semithick, chocolate213940, dashed]
table {%
-2.45 1
51.45 1
};
\addlegendentry{$[E(\{(A,B_\star)\})]_2$}
\end{axis}

\end{tikzpicture}
\caption{Adaptive safe bound.}
\label{fig:adaptive-bound}
\end{subfigure}
\par\vspace{7mm}
\begin{subfigure}[t]{0.47\textwidth}
\centering
\begin{tikzpicture}

\definecolor{chocolate213940}{RGB}{213,94,0}
\definecolor{darkcyan0114178}{RGB}{0,114,178}
\definecolor{darkcyan0158115}{RGB}{0,158,115}
\definecolor{darkgrey176}{RGB}{176,176,176}

\begin{axis}[
    width=0.85\linewidth,
height=0.45\linewidth,
legend cell align={left},
legend columns=2,
legend style={
  fill opacity=0.8,
  draw opacity=1,
  text opacity=1,
        at={(0.5,1.7)},
  anchor=north,
  draw=none
},
tick align=outside,
tick pos=left,
width=\linewidth,
x grid style={darkgrey176},
xlabel={$t$},
xmajorgrids,
xmin=-2.45, xmax=51.45,
xminorgrids,
xtick style={color=black},
y grid style={darkgrey176},
ylabel={$B_{1}$},
ymajorgrids,
ymin=0.03, ymax=0.47,
yminorgrids,
ytick style={color=black}
]
\path [draw=darkcyan0114178, fill=darkcyan0114178, opacity=0.15]
(axis cs:0,0.45)
--(axis cs:0,0.05)
--(axis cs:1,0.05)
--(axis cs:2,0.05)
--(axis cs:3,0.120987052456844)
--(axis cs:4,0.120987052456844)
--(axis cs:5,0.120987052456844)
--(axis cs:6,0.120987052456844)
--(axis cs:7,0.182464431209643)
--(axis cs:8,0.182464431209643)
--(axis cs:9,0.182464431209643)
--(axis cs:10,0.182464431209643)
--(axis cs:11,0.182464431209643)
--(axis cs:12,0.182464431209643)
--(axis cs:13,0.182464431209643)
--(axis cs:14,0.182464431209643)
--(axis cs:15,0.182464431209643)
--(axis cs:16,0.182464431209643)
--(axis cs:17,0.182464431209643)
--(axis cs:18,0.182464431209643)
--(axis cs:19,0.182464431209643)
--(axis cs:20,0.182464431209643)
--(axis cs:21,0.182464431209643)
--(axis cs:22,0.182464431209643)
--(axis cs:23,0.182464431209643)
--(axis cs:24,0.182464431209643)
--(axis cs:25,0.182464431209643)
--(axis cs:26,0.182464431209643)
--(axis cs:27,0.182464431209643)
--(axis cs:28,0.182464431209643)
--(axis cs:29,0.182464431209643)
--(axis cs:30,0.182464431209643)
--(axis cs:31,0.182464431209643)
--(axis cs:32,0.182464431209643)
--(axis cs:33,0.182464431209643)
--(axis cs:34,0.182464431209643)
--(axis cs:35,0.182464431209643)
--(axis cs:36,0.182464431209643)
--(axis cs:37,0.182464431209643)
--(axis cs:38,0.182464431209643)
--(axis cs:39,0.182464431209643)
--(axis cs:40,0.182464431209643)
--(axis cs:41,0.182464431209643)
--(axis cs:42,0.182464431209643)
--(axis cs:43,0.182464431209643)
--(axis cs:44,0.182464431209643)
--(axis cs:45,0.182464431209643)
--(axis cs:46,0.182464431209643)
--(axis cs:47,0.182464431209643)
--(axis cs:48,0.182464431209643)
--(axis cs:49,0.182464431209643)
--(axis cs:49,0.206354965220452)
--(axis cs:49,0.206354965220452)
--(axis cs:48,0.206354965220452)
--(axis cs:47,0.206354965220452)
--(axis cs:46,0.206354965220452)
--(axis cs:45,0.206354965220452)
--(axis cs:44,0.206354965220452)
--(axis cs:43,0.206354965220452)
--(axis cs:42,0.206354965220452)
--(axis cs:41,0.206354965220452)
--(axis cs:40,0.206354965220452)
--(axis cs:39,0.206354965220452)
--(axis cs:38,0.206354965220452)
--(axis cs:37,0.206354965220452)
--(axis cs:36,0.206354965220452)
--(axis cs:35,0.206354965220452)
--(axis cs:34,0.206354965220452)
--(axis cs:33,0.206354965220452)
--(axis cs:32,0.206354965220452)
--(axis cs:31,0.206354965220452)
--(axis cs:30,0.206354965220452)
--(axis cs:29,0.206354965220452)
--(axis cs:28,0.206354965220452)
--(axis cs:27,0.206354965220452)
--(axis cs:26,0.249069517158496)
--(axis cs:25,0.249069517158496)
--(axis cs:24,0.249069517158496)
--(axis cs:23,0.249069517158496)
--(axis cs:22,0.249069517158496)
--(axis cs:21,0.249069517158496)
--(axis cs:20,0.249069517158496)
--(axis cs:19,0.249069517158496)
--(axis cs:18,0.249069517158496)
--(axis cs:17,0.249069517158496)
--(axis cs:16,0.249069517158496)
--(axis cs:15,0.249069517158496)
--(axis cs:14,0.249069517158496)
--(axis cs:13,0.274054752692343)
--(axis cs:12,0.274054752692343)
--(axis cs:11,0.274054752692343)
--(axis cs:10,0.274054752692343)
--(axis cs:9,0.31895068395279)
--(axis cs:8,0.31895068395279)
--(axis cs:7,0.31895068395279)
--(axis cs:6,0.31895068395279)
--(axis cs:5,0.336556131618921)
--(axis cs:4,0.336556131618921)
--(axis cs:3,0.336556131618921)
--(axis cs:2,0.336556131618921)
--(axis cs:1,0.45)
--(axis cs:0,0.45)
--cycle;
\addlegendimage{area legend, draw=darkcyan0114178, fill=darkcyan0114178, opacity=0.15}
\addlegendentry{Interval $\bigl[[\underline B_t]_1,\,[\overline B_t]_1\bigr]$}

\addplot [darkcyan0114178, mark=*, mark size=1, mark options={solid}, forget plot]
table {%
0 0.05
1 0.05
2 0.05
3 0.120987052456844
4 0.120987052456844
5 0.120987052456844
6 0.120987052456844
7 0.182464431209643
8 0.182464431209643
9 0.182464431209643
10 0.182464431209643
11 0.182464431209643
12 0.182464431209643
13 0.182464431209643
14 0.182464431209643
15 0.182464431209643
16 0.182464431209643
17 0.182464431209643
18 0.182464431209643
19 0.182464431209643
20 0.182464431209643
21 0.182464431209643
22 0.182464431209643
23 0.182464431209643
24 0.182464431209643
25 0.182464431209643
26 0.182464431209643
27 0.182464431209643
28 0.182464431209643
29 0.182464431209643
30 0.182464431209643
31 0.182464431209643
32 0.182464431209643
33 0.182464431209643
34 0.182464431209643
35 0.182464431209643
36 0.182464431209643
37 0.182464431209643
38 0.182464431209643
39 0.182464431209643
40 0.182464431209643
41 0.182464431209643
42 0.182464431209643
43 0.182464431209643
44 0.182464431209643
45 0.182464431209643
46 0.182464431209643
47 0.182464431209643
48 0.182464431209643
49 0.182464431209643
};
\addplot [darkcyan0114178, mark=*, mark size=1, mark options={solid}, forget plot]
table {%
0 0.45
1 0.45
2 0.336556131618921
3 0.336556131618921
4 0.336556131618921
5 0.336556131618921
6 0.31895068395279
7 0.31895068395279
8 0.31895068395279
9 0.31895068395279
10 0.274054752692343
11 0.274054752692343
12 0.274054752692343
13 0.274054752692343
14 0.249069517158496
15 0.249069517158496
16 0.249069517158496
17 0.249069517158496
18 0.249069517158496
19 0.249069517158496
20 0.249069517158496
21 0.249069517158496
22 0.249069517158496
23 0.249069517158496
24 0.249069517158496
25 0.249069517158496
26 0.249069517158496
27 0.206354965220452
28 0.206354965220452
29 0.206354965220452
30 0.206354965220452
31 0.206354965220452
32 0.206354965220452
33 0.206354965220452
34 0.206354965220452
35 0.206354965220452
36 0.206354965220452
37 0.206354965220452
38 0.206354965220452
39 0.206354965220452
40 0.206354965220452
41 0.206354965220452
42 0.206354965220452
43 0.206354965220452
44 0.206354965220452
45 0.206354965220452
46 0.206354965220452
47 0.206354965220452
48 0.206354965220452
49 0.206354965220452
};
\addplot [semithick, darkcyan0158115, const plot mark left, mark=*, mark size=1.5, mark options={solid}]
table {%
0 0.45
1 0.05
2 0.336556131618921
3 0.120987052456844
4 0.336556131618921
5 0.120987052456844
6 0.31895068395279
7 0.182464431209643
8 0.31895068395279
9 0.182464431209643
10 0.274054752692343
11 0.182464431209643
12 0.274054752692343
13 0.182464431209643
14 0.249069517158496
15 0.249069517158496
16 0.182464431209643
17 0.249069517158496
18 0.182464431209643
19 0.249069517158496
20 0.182464431209643
21 0.249069517158496
22 0.182464431209643
23 0.249069517158496
24 0.249069517158496
25 0.182464431209643
26 0.249069517158496
27 0.206354965220452
28 0.206354965220452
29 0.206354965220452
30 0.206354965220452
31 0.206354965220452
32 0.206354965220452
33 0.206354965220452
34 0.206354965220452
35 0.206354965220452
36 0.206354965220452
37 0.206354965220452
38 0.206354965220452
39 0.206354965220452
40 0.206354965220452
41 0.206354965220452
42 0.206354965220452
43 0.206354965220452
44 0.206354965220452
45 0.206354965220452
46 0.206354965220452
47 0.206354965220452
48 0.206354965220452
49 0.206354965220452
};
\addlegendentry{$[\widehat B_t]_1$}
\addplot [semithick, chocolate213940, dashed]
table {%
-2.45 0.2
51.45 0.2
};
\addlegendentry{$[B_\star]_1$}
\end{axis}

\end{tikzpicture}
\caption{Consistency interval for $B_1$.}
\label{fig:interval-B1}
\end{subfigure}\hfill
\begin{subfigure}[t]{0.47\textwidth}
\centering
\begin{tikzpicture}

\definecolor{chocolate213940}{RGB}{213,94,0}
\definecolor{darkcyan0114178}{RGB}{0,114,178}
\definecolor{darkcyan0158115}{RGB}{0,158,115}
\definecolor{darkgrey176}{RGB}{176,176,176}

\begin{axis}[
    width=0.85\linewidth,
height=0.45\linewidth,
legend cell align={left},
legend columns=2,
legend style={
  fill opacity=0.8,
  draw opacity=1,
  text opacity=1,
        at={(0.5,1.7)},
  anchor=north,
  draw=none
},
tick align=outside,
tick pos=left,
width=\linewidth,
x grid style={darkgrey176},
xlabel={$t$},
xmajorgrids,
xmin=-2.45, xmax=51.45,
xminorgrids,
xtick style={color=black},
y grid style={darkgrey176},
ylabel={$B_{2}$},
ymajorgrids,
ymin=-0.845, ymax=0.145,
yminorgrids,
ytick style={color=black}
]
\path [draw=darkcyan0114178, fill=darkcyan0114178, opacity=0.15]
(axis cs:0,0.1)
--(axis cs:0,-0.8)
--(axis cs:1,-0.622349988242274)
--(axis cs:2,-0.622349988242274)
--(axis cs:3,-0.470673715579631)
--(axis cs:4,-0.470673715579631)
--(axis cs:5,-0.470673715579631)
--(axis cs:6,-0.470673715579631)
--(axis cs:7,-0.441881062012141)
--(axis cs:8,-0.441881062012141)
--(axis cs:9,-0.441881062012141)
--(axis cs:10,-0.441881062012141)
--(axis cs:11,-0.441881062012141)
--(axis cs:12,-0.441881062012141)
--(axis cs:13,-0.441881062012141)
--(axis cs:14,-0.441881062012141)
--(axis cs:15,-0.441881062012141)
--(axis cs:16,-0.441881062012141)
--(axis cs:17,-0.441881062012141)
--(axis cs:18,-0.441881062012141)
--(axis cs:19,-0.441881062012141)
--(axis cs:20,-0.441881062012141)
--(axis cs:21,-0.441881062012141)
--(axis cs:22,-0.441881062012141)
--(axis cs:23,-0.441881062012141)
--(axis cs:24,-0.441881062012141)
--(axis cs:25,-0.441881062012141)
--(axis cs:26,-0.441881062012141)
--(axis cs:27,-0.441881062012141)
--(axis cs:28,-0.441881062012141)
--(axis cs:29,-0.441881062012141)
--(axis cs:30,-0.441881062012141)
--(axis cs:31,-0.441881062012141)
--(axis cs:32,-0.441881062012141)
--(axis cs:33,-0.441881062012141)
--(axis cs:34,-0.441881062012141)
--(axis cs:35,-0.441881062012141)
--(axis cs:36,-0.441881062012141)
--(axis cs:37,-0.441881062012141)
--(axis cs:38,-0.441881062012141)
--(axis cs:39,-0.441881062012141)
--(axis cs:40,-0.441881062012141)
--(axis cs:41,-0.441881062012141)
--(axis cs:42,-0.441881062012141)
--(axis cs:43,-0.441881062012141)
--(axis cs:44,-0.441881062012141)
--(axis cs:45,-0.441881062012141)
--(axis cs:46,-0.441881062012141)
--(axis cs:47,-0.441881062012141)
--(axis cs:48,-0.441881062012141)
--(axis cs:49,-0.441881062012141)
--(axis cs:49,-0.376657213758312)
--(axis cs:49,-0.376657213758312)
--(axis cs:48,-0.376657213758312)
--(axis cs:47,-0.376657213758312)
--(axis cs:46,-0.376657213758312)
--(axis cs:45,-0.376657213758312)
--(axis cs:44,-0.376657213758312)
--(axis cs:43,-0.376657213758312)
--(axis cs:42,-0.376657213758312)
--(axis cs:41,-0.376657213758312)
--(axis cs:40,-0.376657213758312)
--(axis cs:39,-0.376657213758312)
--(axis cs:38,-0.376657213758312)
--(axis cs:37,-0.355546920959245)
--(axis cs:36,-0.355546920959245)
--(axis cs:35,-0.355546920959245)
--(axis cs:34,-0.355546920959245)
--(axis cs:33,-0.355546920959245)
--(axis cs:32,-0.355546920959245)
--(axis cs:31,-0.355546920959245)
--(axis cs:30,-0.355546920959245)
--(axis cs:29,-0.355546920959245)
--(axis cs:28,-0.355546920959245)
--(axis cs:27,-0.355546920959245)
--(axis cs:26,-0.355546920959245)
--(axis cs:25,-0.227656917924085)
--(axis cs:24,-0.227656917924085)
--(axis cs:23,-0.227656917924085)
--(axis cs:22,-0.227656917924085)
--(axis cs:21,-0.227656917924085)
--(axis cs:20,-0.227656917924085)
--(axis cs:19,-0.227656917924085)
--(axis cs:18,-0.174551775689776)
--(axis cs:17,-0.174551775689776)
--(axis cs:16,-0.174551775689776)
--(axis cs:15,-0.174551775689776)
--(axis cs:14,-0.174551775689776)
--(axis cs:13,-0.174551775689776)
--(axis cs:12,-0.174551775689776)
--(axis cs:11,-0.174551775689776)
--(axis cs:10,-0.174551775689776)
--(axis cs:9,-0.154990036543071)
--(axis cs:8,-0.154990036543071)
--(axis cs:7,-0.154990036543071)
--(axis cs:6,-0.154990036543071)
--(axis cs:5,-0.154990036543071)
--(axis cs:4,-0.154990036543071)
--(axis cs:3,-0.132152264345059)
--(axis cs:2,-0.132152264345059)
--(axis cs:1,0.1)
--(axis cs:0,0.1)
--cycle;
\addlegendimage{area legend, draw=darkcyan0114178, fill=darkcyan0114178, opacity=0.15}
\addlegendentry{Interval $\bigl[[\underline B_t]_2,\,[\overline B_t]_2\bigr]$}

\addplot [darkcyan0114178, mark=*, mark size=1, mark options={solid}, forget plot]
table {%
0 -0.8
1 -0.622349988242274
2 -0.622349988242274
3 -0.470673715579631
4 -0.470673715579631
5 -0.470673715579631
6 -0.470673715579631
7 -0.441881062012141
8 -0.441881062012141
9 -0.441881062012141
10 -0.441881062012141
11 -0.441881062012141
12 -0.441881062012141
13 -0.441881062012141
14 -0.441881062012141
15 -0.441881062012141
16 -0.441881062012141
17 -0.441881062012141
18 -0.441881062012141
19 -0.441881062012141
20 -0.441881062012141
21 -0.441881062012141
22 -0.441881062012141
23 -0.441881062012141
24 -0.441881062012141
25 -0.441881062012141
26 -0.441881062012141
27 -0.441881062012141
28 -0.441881062012141
29 -0.441881062012141
30 -0.441881062012141
31 -0.441881062012141
32 -0.441881062012141
33 -0.441881062012141
34 -0.441881062012141
35 -0.441881062012141
36 -0.441881062012141
37 -0.441881062012141
38 -0.441881062012141
39 -0.441881062012141
40 -0.441881062012141
41 -0.441881062012141
42 -0.441881062012141
43 -0.441881062012141
44 -0.441881062012141
45 -0.441881062012141
46 -0.441881062012141
47 -0.441881062012141
48 -0.441881062012141
49 -0.441881062012141
};
\addplot [darkcyan0114178, mark=*, mark size=1, mark options={solid}, forget plot]
table {%
0 0.1
1 0.1
2 -0.132152264345059
3 -0.132152264345059
4 -0.154990036543071
5 -0.154990036543071
6 -0.154990036543071
7 -0.154990036543071
8 -0.154990036543071
9 -0.154990036543071
10 -0.174551775689776
11 -0.174551775689776
12 -0.174551775689776
13 -0.174551775689776
14 -0.174551775689776
15 -0.174551775689776
16 -0.174551775689776
17 -0.174551775689776
18 -0.174551775689776
19 -0.227656917924085
20 -0.227656917924085
21 -0.227656917924085
22 -0.227656917924085
23 -0.227656917924085
24 -0.227656917924085
25 -0.227656917924085
26 -0.355546920959245
27 -0.355546920959245
28 -0.355546920959245
29 -0.355546920959245
30 -0.355546920959245
31 -0.355546920959245
32 -0.355546920959245
33 -0.355546920959245
34 -0.355546920959245
35 -0.355546920959245
36 -0.355546920959245
37 -0.355546920959245
38 -0.376657213758312
39 -0.376657213758312
40 -0.376657213758312
41 -0.376657213758312
42 -0.376657213758312
43 -0.376657213758312
44 -0.376657213758312
45 -0.376657213758312
46 -0.376657213758312
47 -0.376657213758312
48 -0.376657213758312
49 -0.376657213758312
};
\addplot [semithick, darkcyan0158115, const plot mark left, mark=*, mark size=1.5, mark options={solid}]
table {%
0 0.1
1 -0.622349988242274
2 -0.132152264345059
3 -0.470673715579631
4 -0.154990036543071
5 -0.470673715579631
6 -0.154990036543071
7 -0.441881062012141
8 -0.154990036543071
9 -0.441881062012141
10 -0.174551775689776
11 -0.441881062012141
12 -0.174551775689776
13 -0.441881062012141
14 -0.174551775689776
15 -0.174551775689776
16 -0.441881062012141
17 -0.174551775689776
18 -0.441881062012141
19 -0.227656917924085
20 -0.441881062012141
21 -0.227656917924085
22 -0.441881062012141
23 -0.227656917924085
24 -0.227656917924085
25 -0.441881062012141
26 -0.355546920959245
27 -0.355546920959245
28 -0.355546920959245
29 -0.355546920959245
30 -0.355546920959245
31 -0.355546920959245
32 -0.355546920959245
33 -0.355546920959245
34 -0.355546920959245
35 -0.355546920959245
36 -0.355546920959245
37 -0.355546920959245
38 -0.376657213758312
39 -0.376657213758312
40 -0.376657213758312
41 -0.376657213758312
42 -0.376657213758312
43 -0.376657213758312
44 -0.376657213758312
45 -0.376657213758312
46 -0.376657213758312
47 -0.376657213758312
48 -0.376657213758312
49 -0.376657213758312
};
\addlegendentry{$[\widehat B_t]_2$}
\addplot [semithick, chocolate213940, dashed]
table {%
-2.45 -0.4
51.45 -0.4
};
\addlegendentry{$[B_\star]_2$}
\end{axis}

\end{tikzpicture}
\caption{Consistency interval for $B_2$.}
\label{fig:interval-B2}
\end{subfigure}
\caption{Closed-loop simulation with a fixed unknown
$B_\star\in[\underline B_0,\overline B_0]$. The top panels compare the
adaptive CE controller with the baseline that knows $B_\star$, while
the remaining panels show the adaptive input bound and evolution of the consistency set.}
\label{fig:interval-B-simulation}
\end{figure}
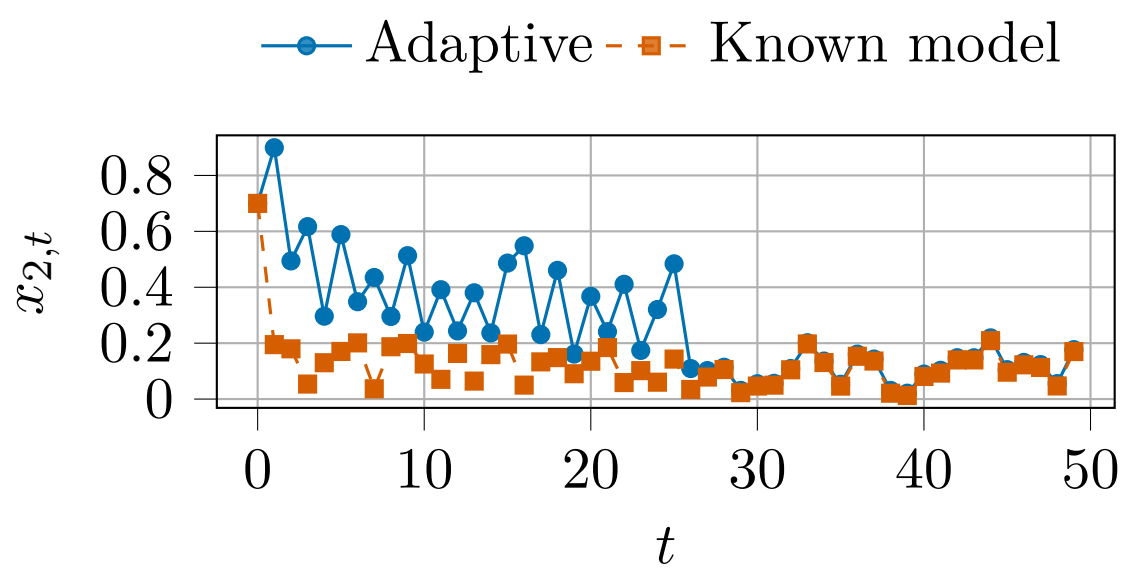
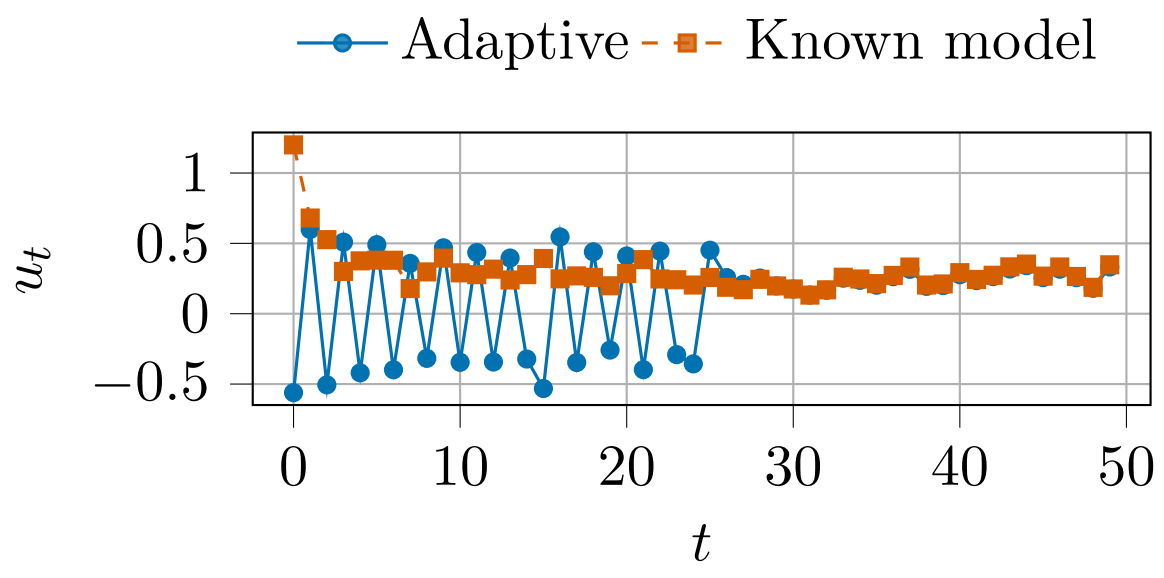
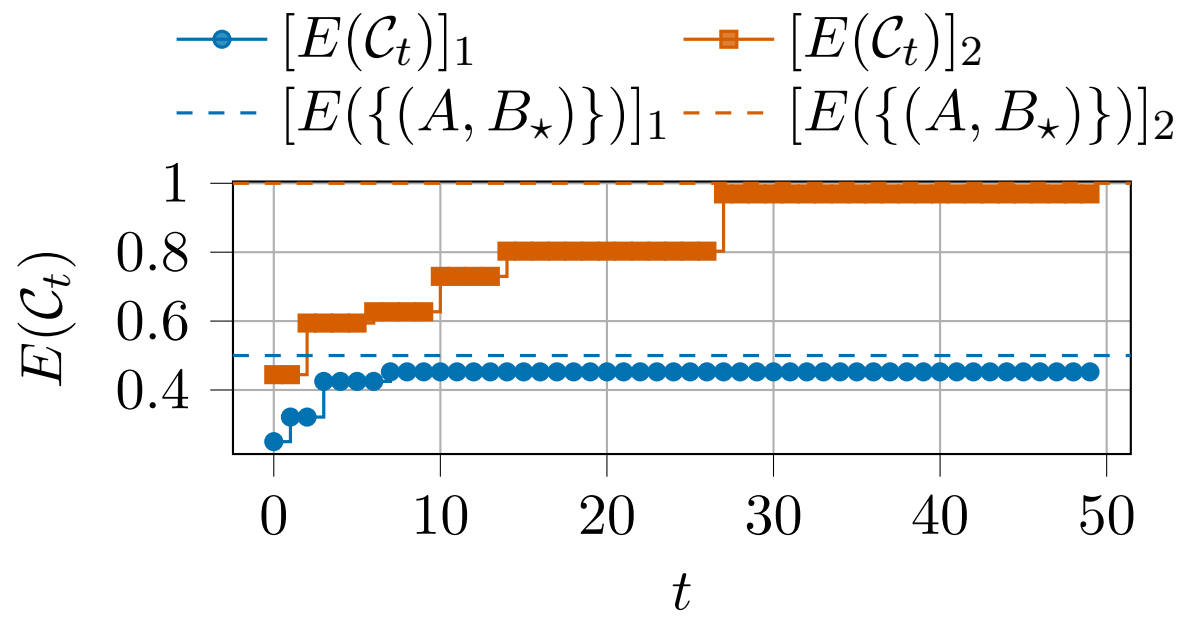
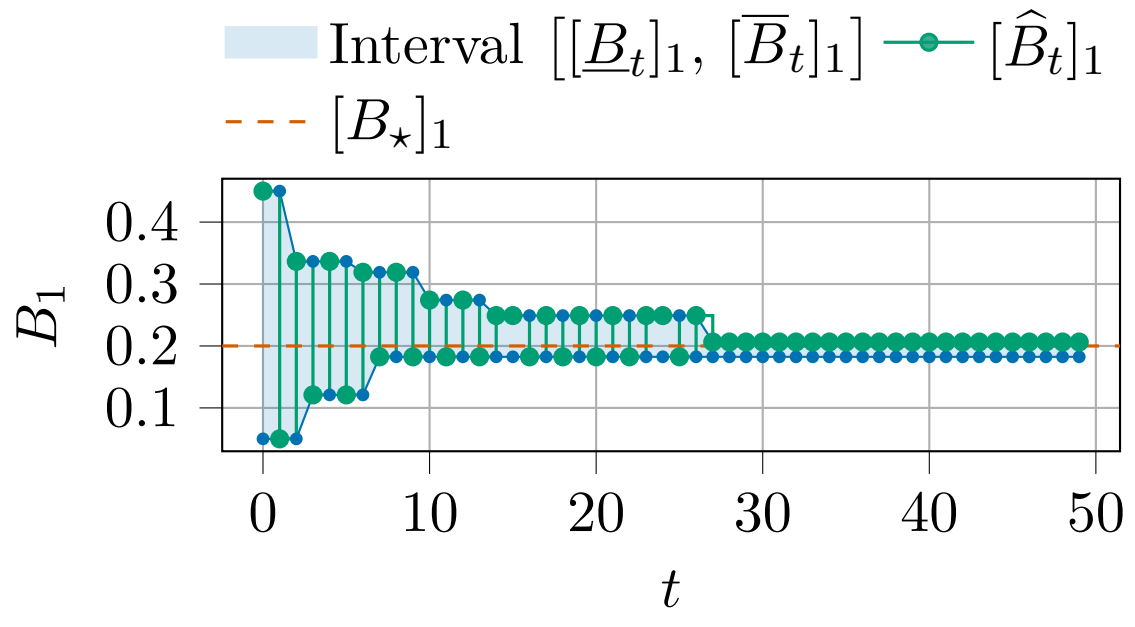
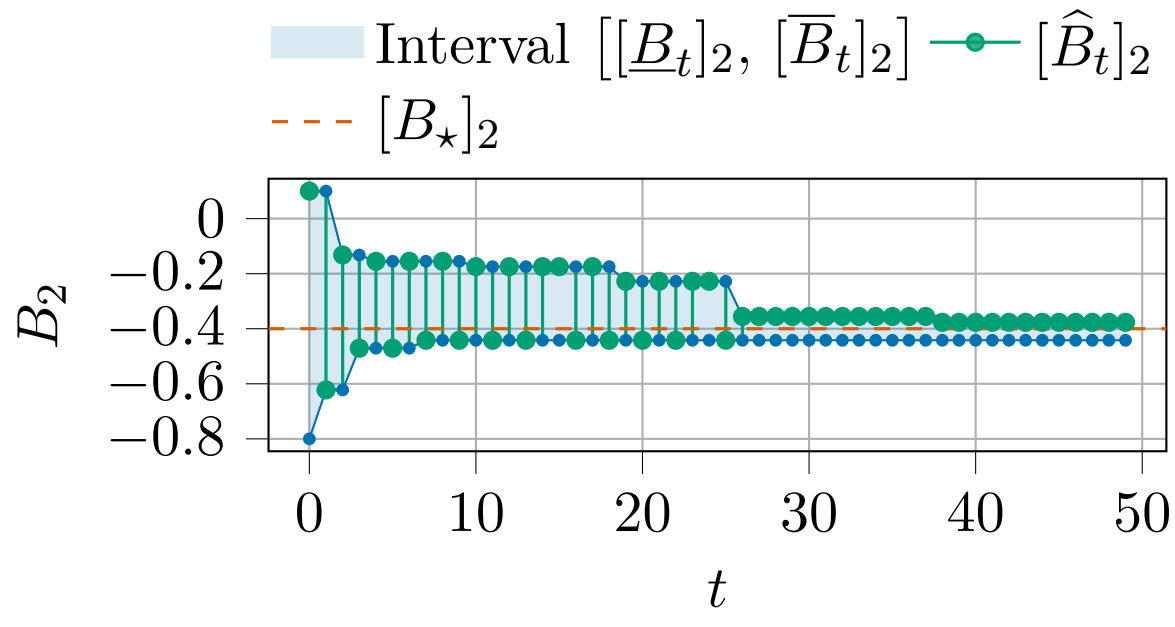

Figure~\ref{fig:interval-B-simulation} shows the resulting trajectories.
The shaded regions in Figures~\ref{fig:interval-B1}
and~\ref{fig:interval-B2} show the componentwise intervals
$[[\underline B_t]_i,[\overline B_t]_i]$, $i\in\{1,2\}$. As observations
become available, $\mathcal C_t$ contracts while continuing to contain the
fixed true pair $(A,B_\star)$. The CE estimate, input bound, and controller
adapt to the refined model set, although $\mathcal C_t$ need not collapse
to a singleton.

The solid curves in Figure~\ref{fig:adaptive-bound} show the adaptive
components of $E(\mathcal C_t)$, whereas the dashed lines show the
known-model components of $E(\{(A,B_\star)\})$. The adaptive bound increases
toward, but need not reach, its known-model value because $\mathcal C_t$
may remain non-singleton. Likewise, the CE
estimate need not converge to $B_\star$. The correct feedback direction is
sufficient, while the remaining gap preserves positivity for all data-consistent
models.

\section{Conclusion}\label{sec:conclusion}
We studied minimax adaptive control of positive systems with fixed but
unknown dynamics and nonnegative disturbances. Using a positivity-induced
consistency set and a historical-data vector, we obtained an exact
reformulation as a minimax dynamic game. Together with the physical state,
these quantities form a sufficient information state for dynamic programming.
The consistency set eliminates models incompatible with observed state
transitions. For finite model sets, we derived an explicit solution to the
Bellman inequality for a certainty-equivalence controller with a
consistency-set-dependent input bound. This bound ensures safe control
throughout adaptation by preserving positivity for every remaining model.
It may expand as incompatible models are eliminated, allowing less cautious
control. We also identified continuum model classes admitting tractable,
fixed-dimensional hyperstate updates, supporting scalability with respect
to the number of candidate model pairs in the uncertainty set.
This work suggests several directions for future research, including:
\begin{itemize}
	\item[\textit{(i)}] Extend the Bellman analysis and performance guarantees
	to continuum uncertainty classes, including the interval classes in
	Examples~\ref{example1}--\ref{example3}.
	
	\item[\textit{(ii)}] Identify broader continuum model classes that admit
	fixed-dimensional recursive hyperstate representations and exploit this
	structure to derive explicit adaptive dual control policies with certified
	performance guarantees.
\end{itemize}

\section*{Acknowledgments}

The authors would like to thank Anders Rantzer for valuable discussions that helped improve the paper.
The authors are affiliated with the ELLIIT Strategic Research Area
at Lund University. This project received funding from the European Research
Council (ERC) under grant agreement No.~834142 (ScalableControl) and No. 101199738 (DualControl), and was
partially supported by the Wallenberg AI, Autonomous Systems and Software
Program (WASP), funded by the Knut and Alice Wallenberg Foundation.

\appendix
\section{Proofs}

 \subsection{Proof of Theorem~\ref{thm:minimax-dp}}\label{proofappenA}
The proof follows \cite[Appendix C.1]{Bencherki2026}, with the information
state replaced by $(x,\mathcal C,Z)$.
Recall the hyperstate update in~\eqref{eq:dp-hyperstate-update}
where $\mathcal C_+\neq\varnothing$. Recall also the Bellman operator~\eqref{eq:Bellman-operator} and the value iteration sequence
$\{V_k\}_{k\in\mathbb N}$ in
\eqref{eq:V0}--\eqref{eq:value-iteration}.

\emph{Proof of \textit{(i)}--\textit{(ii)}.}
We start by showing that the sequence
$V_0,V_1,V_2,\ldots$ is monotonically nondecreasing. By definition,
\begin{equation*}
V_1(x,\mathcal C,Z)
=
\min_{u\in\mathcal U(x,\mathcal C)}
\max_{\substack{v,
\mathcal C_+\neq\varnothing}}
\Bigl\{
s^\top x+r^\top u
-
\min_{(A,B)\in\mathcal C_+}
\gamma^\top
[\,I\ {-}A\ {-}B\,]Z_+
\Bigr\}.
\end{equation*}

Choose a model $(A,B)\in\mathcal C$ attaining the minimum
that defines $V_0$, so that
\[
\gamma^\top[\,I\ {-}A\ {-}B\,]Z
=
\min_{(A',B')\in\mathcal C}
\gamma^\top[\,I\ {-}A'\ {-}B'\,]Z
=
-V_0(x,\mathcal C,Z).
\]
Such a model exists because $\mathcal C$ is nonempty and compact.
For this fixed model and every $u\in\mathcal U(x,\mathcal C)$,
the choice $v=Ax+Bu$ is admissible and ensures
$(A,B)\in\mathcal C_+$ and $\mathcal C_+\neq\varnothing$.
Moreover, this choice of $v$ gives
\begin{equation*}
\gamma^\top[\,I\ {-}A\ {-}B\,]Z_+
=
\gamma^\top[\,I\ {-}A\ {-}B\,]Z
+\gamma^\top(v-Ax-Bu)
=
\gamma^\top[\,I\ {-}A\ {-}B\,]Z.
\end{equation*}
Since the selected model belongs to $\mathcal C_+$,
the minimum over $\mathcal C_+$ cannot exceed its value
at this model. Therefore,
\begin{multline*}
-\min_{(A',B')\in\mathcal C_+}
\gamma^\top[\,I\ {-}A'\ {-}B'\,]Z_+
\geq
-\gamma^\top[\,I\ {-}A\ {-}B\,]Z_+
=
-\gamma^\top[\,I\ {-}A\ {-}B\,]Z
\\=
-\min_{(A',B')\in\mathcal C}
\gamma^\top[\,I\ {-}A'\ {-}B'\,]Z
=
V_0(x,\mathcal C,Z).
\end{multline*}
Thus, for every admissible $u$, there is an admissible $v$
such that $V_0(v,\mathcal C_+,Z_+)\geq V_0(x,\mathcal C,Z)$.
Consequently, maximizing $V_0(v,\mathcal C_+,Z_+)$ over all
admissible $v$ gives a value at least $V_0(x,\mathcal C,Z)$.
It follows that
\begin{multline*}
V_1(x,\mathcal C,Z)
\geq
\min_{u\in\mathcal U(x,\mathcal C)}
\Bigl\{
s^\top x+r^\top u
-
\min_{(A,B)\in{\mathcal C}}
\gamma^\top[\,I\ {-}A\ {-}B\,]Z
\Bigr\}
\geq
-\min_{(A,B)\in\mathcal C}
\gamma^\top[\,I\ {-}A\ {-}B\,]Z
\\=
V_0(x,\mathcal C,Z),
\end{multline*}
where the second inequality (positivity of the stage cost) follows from
Assumption~\ref{ass:positive-stage-cost}, since
$s^\top x+r^\top u\geq0$ for every admissible $u$.
Thus, $V_1\geq V_0$.
Since $\mathcal F$ is order preserving in $V$, it
follows by induction that
$V_{k+1}=\mathcal F V_k\geq\mathcal F V_{k-1}=V_k$ for $k\geq1$.
We conclude that the sequence
$V_0,V_1,V_2,\ldots$ is monotonically nondecreasing.

Fix $N\geq0$. Consider the truncated original minimax problem
\begin{gather*}
\inf_\mu
\sup_{\substack{
(A,B)\in\mathcal M\\
w\geq0}}
\sum_{t=0}^{N}
\bigl(
s^\top x_t+r^\top u_t-\gamma^\top w_t
\bigr)\\
\mathrm{s.t.}\quad
x_{t+1}=Ax_t+Bu_t+w_t,
\qquad
x_0\in\mathbb R_+^n.
\end{gather*}
The value of Problem~\ref{pb1} is bounded below by the value of this
truncated problem. The truncated value grows monotonically with $N$,
and the value of Problem~\ref{pb1} is obtained in the limit as
$N\to\infty$.
Introducing $v_t=x_{t+1}$ and using the consistency set and history
vector renders the truncated problem equivalent to
\begin{gather*}
	\inf_\eta
	\sup_{\substack{v\\ \mathcal C_{N+1}\neq\varnothing}}
	\Bigl\{
	\sum_{t=0}^{N}
	\left(s^\top x_t+r^\top u_t\right)
	-
	\min_{(A,B)\in\mathcal C_{N+1}}
	\gamma^\top[\,I\ {-}A\ {-}B\,]Z_{N+1}
	\Bigr\}\\
	\begin{aligned}
		\mathrm{s.t.}\quad x_{t+1}
		&=v_t,\\
		\mathcal C_{t+1}
		&=\bigl\{
		(A,B)\in\mathcal C_t:\,
		v_t\geq Ax_t+Bu_t
		\bigr\},
		\quad \mathcal C_0=\mathcal M,\\
		Z_{t+1}
		&=Z_t+
		\begin{bmatrix}
			v_t^\top & x_t^\top & u_t^\top
		\end{bmatrix}^{\top},
		\quad Z_0=0,\\
		u_t
		&=\eta(x_t,\mathcal C_t,Z_t),
		\quad u_t\in\mathcal U(x_t,\mathcal C_t).
	\end{aligned}
\end{gather*}
Standard minimax dynamic programming with terminal cost
\eqref{eq:V0} and recursion \eqref{eq:value-iteration} shows
that the truncated problem has value
$V_{N+1}(x_0,\mathcal M,0)$. Thus, Problem~\ref{pb1} has
finite value if and only if
$\{V_k(x_0,\mathcal M,0)\}_{k\geq0}$ is bounded above,
and $J_\star(x_0)=V_\star(x_0,\mathcal M,0)$.
For every reachable $(x,\mathcal C,Z)$,
$V_k(x,\mathcal C,Z)\leq V_k(x,\mathcal C,0)
\leq V_k(x,\mathcal M,0)$,
since accumulated disturbance penalties are nonnegative and
the admissible input set can only expand as models are ruled out.
Thus, boundedness at $(x,\mathcal M,0)$ for every $x\geq0$
also ensures a finite limit at every reachable hyperstate. If Problem~\ref{pb2} has finite value, then
$V_k(x_0,\mathcal M,0)$ is bounded above by that value,
so also
$V_\star(x_0,\mathcal M,0)
=\lim_{k\to\infty}V_k(x_0,\mathcal M,0)$
is finite.
The minimizing argument $\eta^\star(x,\mathcal C,Z)$ of
$\mathcal F V_\star$ is then optimal for Problem~\ref{pb2},
and the induced policy $\mu^\star$ defined by
\eqref{eq:mu-star}--\eqref{eq:Ct-history} is optimal for
Problem~\ref{pb1}. This concludes the proof of
\textit{(i)}--\textit{(ii)}.

\emph{Proof of \textit{(iii)}.}
Suppose that $\overline V\geq V_0$ and
$
\mathcal F_{\overline\eta(x,\mathcal{C},Z)}\overline V
\leq\overline V.
$
Define the sequence $W_0,W_1,W_2,\ldots$ recursively by
\begin{align}
W_0\coloneqq V_0,\quad
W_{k+1}(x,\mathcal C,Z)
\coloneqq
\mathcal F_{\overline\eta(x,\mathcal{C},Z)}
W_k(x,\mathcal C,Z).
\end{align}
By dynamic programming,
\begin{equation*}
W_{N+1}(x_0,\mathcal M,0)
=
\sup_{\substack{
v\\
\mathcal C_{N+1}\neq\varnothing}}
\Bigl\{
\sum_{t=0}^{N}
\left(
s^\top x_t+r^\top u_t
\right)
-
\min_{(A,B)\in\mathcal C_{N+1}}
\gamma^\top
[\,I\ {-}A\ {-}B\,]
Z_{N+1}
\Bigr\},
\end{equation*}
where
$u_t=\overline\eta(x_t,\mathcal C_t,Z_t)$, and the hyperstate is
generated by \eqref{eq:dp-hyperstate-update}. By induction, the definitions of \(V_k\) and \(W_k\) imply that $\overline V\geq W_k\geq V_k$ for all $k\geq0$.
Indeed, if $\overline V \ge W_k\ge V_k$ for some $k$,
\begin{equation*}
V_{k+1}
=
\mathcal F V_k
\leq
\mathcal F W_k
\leq
\mathcal F_{\overline\eta(x,\mathcal{C},Z)}W_k
=
W_{k+1}
\leq
\mathcal F_{\overline\eta(x,\mathcal{C},Z)}\overline V
\leq
\overline V.
\end{equation*}
Consequently,
$
V_\star(x,\mathcal C,Z)
\leq
\lim_{k \to \infty }W_k(x,\mathcal C,Z)
\leq
\overline V(x,\mathcal C,Z).
$
This proves that
$
V_\star(x_0,\mathcal M,0)
\leq J_{\overline\mu}(x_0)
\leq
\overline V(x_0,\mathcal M,0)
$
when the control law $\overline\mu$ is defined by
$
u_t=\overline\eta(x_t,\mathcal C_t,Z_t)
$, with $\mathcal C_t$
and $Z_t$ generated according to \eqref{eq:dp-hyperstate-update}.
In particular, $\overline\mu$ is optimal if
$\overline V(x,\mathcal C,Z)=V_\star(x,\mathcal C,Z)$.
{\flushright\qedblack}

\subsection{Proof of Theorem~\ref{thm2}}
\label{thm2-proof}

We prove Theorem~\ref{thm2} by verifying item~\textit{(iii)} of
Theorem~\ref{thm:minimax-dp}.
We start by showing that $\overbar V\geq V_0$. If
$\mathcal C=\{(A_i,B_i)\}$, so that $|\mathcal C|=1$, then
$p_i^\top x\geq0$ gives
\[
\begin{aligned}
\overbar V(x,\mathcal C,Z)
=
p_i^\top x
-\gamma^\top[\,I\ {-}A_i\ {-}B_i\,]Z
\geq
-\gamma^\top[\,I\ {-}A_i\ {-}B_i\,]Z
=
V_0(x,\mathcal C,Z).
\end{aligned}
\]
Alternatively, if $|\mathcal C|>1$, retaining only the terms
with $i=j$ in the maximum defining $\overbar V$ gives
\begin{equation*}
\overbar V(x,\mathcal C,Z)
\geq
\max_{i:(A_i,B_i)\in\mathcal C}
\Bigl\{
p_{ii}^\top x
-\gamma^\top[\,I\ {-}A_i\ {-}B_i\,]Z
\Bigr\}
\geq V_0(x,\mathcal C,Z),
\end{equation*}
where the last inequality uses $p_{ii}^\top x\geq0$.
Thus, $\overbar V\geq V_0$.

We now show that $\mathcal F\overbar V
\leq\mathcal F_{\overbar\eta}\overbar V
\leq\overbar V$. Recall the policy
\begin{equation*}
u=\overbar\eta(x,\mathcal C,Z)=K_k(\mathcal C)x,
\;\;
k\in\operatorname*{arg\,min}_{i:(A_i,B_i)\in\mathcal C}
\gamma^\top[\,I\ {-}A_i\ {-}B_i\,]Z.
\end{equation*}
Since
$
|u|
\leq
|K_k(\mathcal C)|x
\leq
E(\mathcal C)x,
$
the input belongs to $\mathcal U(x,\mathcal C)$. Moreover, for every
$(A_i,B_i)\in\mathcal C$,
$
A_ix+B_iu
\geq
\bigl(A_i-|B_i|E(\mathcal C)\bigr)x
\geq0.
$
Hence the selected input preserves positivity for every model in
$\mathcal C$. The choice of $k$ also gives
\begin{equation}
\gamma^\top[\,I\ {-}A_k\ {-}B_k\,]Z
\leq
\gamma^\top[\,I\ {-}A_i\ {-}B_i\,]Z,
\quad \forall\,(A_i,B_i)\in\mathcal C.
\label{eq:selected-model-order}
\end{equation}
For a feasible next state $v$, recall the hyperstate update
\begin{equation}
\begin{aligned}
\mathcal C_+
&=
\{(A,B)\in\mathcal C:
v\geq Ax+BK_k(\mathcal C)x\},\\
Z_+
&=
Z+
\begin{bmatrix}
v^\top & x^\top & (K_k(\mathcal C)x)^\top
\end{bmatrix}^{\!\top}.
\end{aligned}
\label{eq:Bellman-updates-substituted}
\end{equation}

Suppose first that $\mathcal C=\{(A_k,B_k)\}$, i.e., $|\mathcal C|=1$. Then
$K_k(\mathcal C)=K_k$ and
$\mathcal C_+=\mathcal C$. Therefore,
\begin{multline}
s^\top x+r^\top K_kx
+\overbar V(v,\mathcal C_+,Z_+)
=
-\gamma^\top[\,I\ {-}A_k\ {-}B_k\,]Z
+s^\top x+r^\top K_kx
+(p_k-\gamma)^\top v
+\gamma^\top(A_k+B_kK_k)x\\
\leq
-\gamma^\top[\,I\ {-}A_k\ {-}B_k\,]Z
+s^\top x+r^\top K_kx
+p_k^\top(A_k+B_kK_k)x
=
-\gamma^\top[\,I\ {-}A_k\ {-}B_k\,]Z
+p_k^\top x
\\=
\overbar V(x,\mathcal C,Z).
\label{eq:proof-singleton}
\end{multline}
The inequality uses $p_k\leq\gamma$ and
$v\geq(A_k+B_kK_k)x\geq0$, whereas the following equality follows from
Assumption~\ref{ass:model-dependent-bounds}\textit{(iii)}.

It remains to consider $|\mathcal C|>1$. Suppose first that the next
observation leaves only one model, so
$\mathcal C_+=\{(A_i,B_i)\}$. The successor value is then
\[
\overbar V(v,\mathcal C_+,Z_+)
=
p_i^\top v-\gamma^\top[\,I\ {-}A_i\ {-}B_i\,]Z_+.
\]
The input producing this observation is still $K_k(\mathcal C)x$.
The model-based gain $K_i$ is used only after this transition.
Using the update of $Z$, we obtain
\begin{multline}
s^\top x+r^\top K_k(\mathcal C)x
+\overbar V(v,\mathcal C_+,Z_+)
=
-\gamma^\top[\,I\ {-}A_i\ {-}B_i\,]Z
+s^\top x+r^\top K_k(\mathcal C)x
+(p_i-\gamma)^\top v
\\+\gamma^\top(A_i+B_iK_k(\mathcal C))x
\leq
-\gamma^\top[\,I\ {-}A_i\ {-}B_i\,]Z
+s^\top x+r^\top K_k(\mathcal C)x
+p_i^\top(A_i+B_iK_k(\mathcal C))x
\\\leq
-\gamma^\top[\,I\ {-}A_i\ {-}B_i\,]Z
+s^\top x+r^\top K_k(\mathcal C)x
+p_{ii}^\top(A_i+B_iK_k(\mathcal C))x
\leq
-\gamma^\top[\,I\ {-}A_i\ {-}B_i\,]Z
\\+p_{ik}^\top x
\leq
V_{ik}(x,Z)
\leq
\overbar V(x,\mathcal C,Z).
\label{eq:proof-singleton-successor}
\end{multline}
The first inequality uses $p_i\leq\gamma$ and
$v\geq(A_i+B_iK_k(\mathcal C))x$. The second uses
$p_i\leq p_{ii}$ and positivity of
$A_i+B_iK_k(\mathcal C)$. The third follows from
\eqref{eq:finite-family-inequality} with $j=i$, and the fourth follows
from \eqref{eq:selected-model-order}. Finally,
$(A_i,B_i)\in\mathcal C_+\subseteq\mathcal C$ and
$(A_k,B_k)\in\mathcal C$. Thus, the pair $(i,k)$ is included
in the maximum defining $\overbar V(x,\mathcal C,Z)$,
which gives the last inequality.

Now suppose that $|\mathcal C_+|>1$. By definition,
\begin{equation}
\overbar V(v,\mathcal C_+,Z_+)
=
\max_{\substack{i,j:
(A_i,B_i),(A_j,B_j)\in\mathcal C_+}}
V_{ij}(v,Z_+).
\label{eq:successor-value}
\end{equation}
Fix any $(A_i,B_i),(A_j,B_j)\in\mathcal C_+$. Their consistency with $v$
implies
\begin{equation}
v\geq
\max_{\ell\in\{i,j\}}
(A_\ell+B_\ell K_k(\mathcal C))x.
\label{eq:pair-feasibility}
\end{equation}
Using the update of $Z$ gives
\begin{gather*}
s^\top x+r^\top K_k(\mathcal C)x
+V_{ij}(v,Z_+)
=
-\tfrac12\gamma^\top
\Bigl(
[\,I\ {-}A_i\ {-}B_i\,]
+[\,I\ {-}A_j\ {-}B_j\,]
\Bigr)Z
+s^\top x+r^\top K_k(\mathcal C)x\\
+(p_{ij}-\gamma)^\top v
+\tfrac12\gamma^\top
\bigl[
A_i+A_j
+(B_i+B_j)K_k(\mathcal C)
\bigr]x.
\end{gather*}
Since $p_{ij}\leq\gamma$, \eqref{eq:pair-feasibility} yields
$
(p_{ij}-\gamma)^\top v
\leq
(p_{ij}-\gamma)^\top
\max_{\ell\in\{i,j\}}
(A_\ell+B_\ell K_k(\mathcal C))x.
$
Furthermore,
\begin{gather*}
2\max_{\ell\in\{i,j\}}
(A_\ell+B_\ell K_k(\mathcal C))x
-\bigl[
A_i+A_j
+(B_i+B_j)K_k(\mathcal C)
\bigr]x
=
\left|
\bigl[
A_i-A_j
+(B_i-B_j)K_k(\mathcal C)
\bigr]x
\right|.
\end{gather*}
Combining these relations gives
\begin{gather*}
s^\top x+r^\top K_k(\mathcal C)x
+V_{ij}(v,Z_+)
\leq s^\top x+r^\top K_k(\mathcal C)x
-\tfrac12\gamma^\top
\Bigl(
[\,I\ {-}A_i\ {-}B_i\,]
+[\,I\ {-}A_j\ {-}B_j\,]
\Bigr)Z\\
+p_{ij}^\top
\max_{\ell\in\{i,j\}}
(A_\ell+B_\ell K_k(\mathcal C))x
-\tfrac12\gamma^\top
\left|
\bigl[
A_i-A_j
+(B_i-B_j)K_k(\mathcal C)
\bigr]x
\right|.
\end{gather*}
If $i\neq j=k$, interchange $i$ and $j$ in
\eqref{eq:finite-family-inequality}; otherwise, apply that condition
directly. The symmetry of $p_{ij}$, the elementwise maximum, and the
absolute value then give
\begin{multline}
s^\top x+r^\top K_k(\mathcal C)x
+V_{ij}(v,Z_+)
\leq
-\tfrac12\gamma^\top
\Bigl(
[\,I\ {-}A_i\ {-}B_i\,]
+[\,I\ {-}A_j\ {-}B_j\,]
\Bigr)Z
+\max\left\{
p_{ik}^\top x,
p_{jk}^\top x
\right\}\\
\leq
\max\left\{
V_{ik}(x,Z),
V_{jk}(x,Z)
\right\}
\leq
\overbar V(x,\mathcal C,Z).
\label{eq:proof-pair-final}
\end{multline}
The second inequality follows from \eqref{eq:selected-model-order}, and
the last follows because $i,j,k$ index models in $\mathcal C$.

Taking the maximum over $i,j$ in \eqref{eq:proof-pair-final} and using
\eqref{eq:successor-value} proves the required bound when
$|\mathcal C_+|>1$. Together with \eqref{eq:proof-singleton} and
\eqref{eq:proof-singleton-successor}, we therefore have
\[
s^\top x+r^\top K_k(\mathcal C)x
+\overbar V(v,\mathcal C_+,Z_+)
\leq
\overbar V(x,\mathcal C,Z)
\]
for every feasible next state $v$. Maximizing over all such $v$ gives
$\mathcal F_{\overbar\eta(x,\mathcal C,Z)}\overbar V\leq\overbar V$.
Since $\overbar\eta(x,\mathcal C,Z)$ is admissible and $\mathcal F$ minimizes over all
admissible inputs,
$
\mathcal F\overbar V
\leq
\mathcal F_{\overbar\eta(x,\mathcal C,Z)}\overbar V
\leq
\overbar V.
$

Finally, $Z_0=0$ and $\mathcal C_0=\mathcal M$ give
$
\overbar V(x_0,\mathcal M,0)
\leq
\max_{i,j\in\{1,\ldots,M\}}p_{ij}^\top x_0,
$
where the inequality follows from the definition of
$\overbar V$ and $p_i\leq p_{ii}$.
Item~\textit{(iii)} of
Theorem~\ref{thm:minimax-dp} now gives
$J_{\overbar\mu}(x_0)
\leq
\max_{i,j\in\{1,\ldots,M\}}p_{ij}^\top x_0$, completing the proof.

\hfill\qedblack

\printbibliography
\end{document}